\documentstyle[amssymb,amsfonts]{amsart}

\newcommand{\Vol}{{\rm Vol}}
\newcommand{\dist}{\operatorname{dist}}
\newcommand{\supp}{\operatorname{supp}}

\def\R{{\hbox{\bf R}}}

\def\Z{{\hbox{\bf Z}}}

\def\be#1{\begin{equation}\label{#1}}
\def\bas{\begin{align*}}
\def\eas{\end{align*}}
\def\bi{\begin{itemize}}
\def\ei{\end{itemize}}

\def\chr{{\hbox{char}}}
\def\Im{{\hbox{Im}}}

\def\eps{\varepsilon}
\newenvironment{proof}{\noindent {\bf Proof} }{\endprf\par}
\def \endprf{\hfill  {\vrule height6pt width6pt depth0pt}\medskip}
\def\emph#1{{\it #1}}
\def\textbf#1{{\bf #1}}

\theoremstyle{plain}
  \newtheorem{theorem}[subsection]{Theorem}

  \newtheorem{proposition}[subsection]{Proposition}
  
  \newtheorem{lemma}[subsection]{Lemma}
  \newtheorem{corollary}[subsection]{Corollary}

\theoremstyle{remark}

\theoremstyle{definition}
  \newtheorem{definition}[subsection]{Definition}

\include{psfig}

\begin{document}

\title{Restriction and Kakeya phenomena for finite fields}

\author{Gerd Mockenhaupt}
\address{Georgia Institute of Technology Atlanta, GA 30332-0160}
\email{gerdm@@math.gatech.edu}

\author{Terence Tao}
\address{Department of Mathematics, UCLA, Los Angeles CA 90095-1555}
\email{tao@@math.ucla.edu}

\subjclass{42B05, 11T24, 52C17}

\begin{abstract} 
The restriction and Kakeya problems in Euclidean space have 
received much attention in the last few decades, and are related to many 
problems in harmonic analysis, PDE, and number theory.  In this paper we 
initiate the study of these problems on finite fields.  In many cases the 
Euclidean arguments carry over easily to the finite setting (and are in fact 
somewhat cleaner), but there are some new phenomena in the finite case which 
deserve closer study. 
\end{abstract}

\maketitle

\centerline{\emph{In memory of Tom Wolff}}

\centerline{}

\section{Introduction}\label{intro-sec}
The purpose of this paper is to begin an investigation of restriction and 
Kakeya phenomena in the space $F^n$, 
where $F$ is a finte field of characteristic $\chr (F)>2$.  
The corresponding phenomena in 
Euclidean space $\R^n$ have been intensively studied (see e.g. Tom Wolff's
survey of the Kakeya phenomenon \cite{wolff:survey} , and 
\cite{stein:problem}, \cite{stein:large}, \cite{borg:stein}, \cite{Bo}
for a discussion of the restriction phenomenon) but it seems the finite 
field case has received far less attention.  This is unfortunate, since 
the finite field case serves as a good model for the Euclidean case in 
that many of the technical difficulties (small angle issues, small 
separation issues, a logarithmically infinite number of scales, and 
Schwartz tails arising from the uncertainty principle) are 
eliminated.  Also these finite field problems are closely related to 
other well-studied problems in number theory and arithmetic or geometric 
combinatorics, and techniques from these areas of mathematics may well be 
useful in attacking these problems.
On the other hand, certain Euclidean space tools are not available in the 
finite field setting (e.g. Taylor series approximations, induction on 
scales, or combinatorial arguments based on the ordering of $\R$). 

In this introduction we only give a brief overview of our results.
Restriction estimates will be discussed in detail in Section \ref{rest-sec} 
onwards, while Kakeya estimates will be discussed in detail in Section 
\ref{kakeya-sec}.  Finally, we show a partial connection between the 
restriction and Kakeya problems in Section \ref{conn-sec}.  Unfortunately 
the connection is not as tight as in the Euclidean case because one does 
not automatically have Taylor series approximations.

Restriction phenomena address the question of estimating the size of 
Fourier transforms of measures such as $fd\sigma$, where $d\sigma$ is 
surface measure on some set $S$ and $f$ is an arbitrary $L^p$ function.  
This problem in Euclidean space was first posed and partially solved 
by Stein in 1967 \cite{feff:thesis}.  
A typical inequality in Euclidean space is the Tomas-Stein inequality
$$ 
\| (f d\sigma)\spcheck \|_{L^{2(n+1)/(n-1)}(\R^n)} \,\leq \, 
C_n \ \| f \|_{L^2(S^{n-1})}
$$
where $f$ is an $L^2$ function on the sphere $S^{n-1}$ and $d\sigma$ 
is surface measure on the sphere.

In this paper we study this problem in the finite field setting for 
various algebraic varieties $S$; there are many sets $S$ of interest, 
but we shall mostly restrict our attention to cones and paraboloids.  
In two and three dimensions we shall be able to get reasonably good 
results, and in particular we can improve upon the standard Tomas-Stein 
restriction theorem for the paraboloid in three dimensions. Moreover,
our results give reason to conjecture that in case $-1$ is not 
a square in $F$ the paraboliods are {\it sharp} $\Lambda(3)$-sets.

The Kakeya problem addresses the problem of the extent to which lines in 
different directions can overlap.  (There are many interesting variants 
of this problem in which lines are replaced by other geometric objects, 
but we will not discuss them here).  We will not prove too many new results 
for these problems, but we present some simplified proofs of various 
Euclidean arguments for the finite field setting, and show why this 
problem is connected to restriction problems and also to more classical 
problems in incidence geometry combinatorics.

This paper is partly of an expository nature; we have tried to make it 
accessible to readers who are not familiar with the restriction and 
Kakeya problems in the Euclidean case.  We hope in particular that 
experts in combinatorics or number theory will be able to shed insight 
on these problems using techniques from those fields.

This work was conducted at UNSW and the University of Milano-Bicocca.  
The second author is a Clay Prize Fellow and is supported by the Packard 
Foundation.

\section{General notation}\label{notation-sec}

If $A$ is a finite set, we use $|A|$ to denote the cardinality of $A$.  
We use $\delta_i$ to denote the Kronecker delta point mass at $i$, thus 
$\delta_i(j) = 1$ if $i=j$ and $\delta_i(j) = 0$ otherwise.

If $p$ is an exponent, we use $p'$ to denote the dual exponent 
$1/p + 1/p' := 1$.

For positive numbers x and y we use $x \lesssim y$ or $x = O(y)$ to 
denote the estimate $x \leq C y$, where $C$ is a constant which depends on 
the dimension $n$ (and possibly on exponents such as $p$, $q$) but not on 
the underlying field $F$.  We also use the notation $x \lessapprox y$ to 
denote the estimate $x \leq C_\eps\, |F|^\eps\, y$ for all $\eps > 0$, 
where $C_\eps$ is a constant which can depend on $\eps$ 
but not on $F$.  This notation is convenient for suppressing logarithmic 
powers of $F$.  Unlike the Euclidean case, we do not have a logarithmically 
infinite number of scales, and so we expect that many of the logarithmic 
losses in this paper can in fact be removed.  On the other hand, many of
our estimates fail to be sharp by an actual power of $|F|$ rather than 
just an epsilon, and so the logarithmic losses are fairly insignificant 
by comparison in these cases.

We shall frequently use the following easy consequence of the Cauchy-Schwarz 
inequality: if $P$, $L$ are two finite sets with some relation $p \sim l$ 
between elements $p \in P$, $l \in L$, then
\be{cauchy}
|\{ (p, l, l') \in P \times L \times L: p \sim l, p \sim l' \}|
\geq \frac{ |\{ (p, l) \in P \times L: p \sim l \}|^2 } {|P|}.
\end{equation} 

We shall be working in the vector space $F^n$, which we endow with counting 
measure $dx$, and the dual space $F^n_*$, 
which we endow with normalized counting measure $d\xi$ which assigns a 
mass of $|F|^{-n}$ to each point, so that $F^n_*$ has total mass $1$.
Thus
$$ \int_{F^n} f(x)\ dx = \sum_{x \in F^n} f(x) \hbox{ and }
\int_{F_*^n} g(\xi)\ d\xi = \frac 1{|F|^n} \sum_{\xi \in F^n_*} g(\xi).$$
The use of integrals for finite sums may appear a bit strange but this is 
a convenient normalization to prevent distracting powers of $|F|$ 
from appearing in the computations.

We fix $e: F \to S^1$ to be a non-principal character of $F$, i.e. a 
multiplicative function from $F$ to the unit circle which is not identically 1. 
For instance, if $F = \Z/p\Z$ for some prime $p$, one can take 
$e(x) := \exp(2\pi i x/p)$.  In practice the exact choice of this character 
will be irrelevant.

For any complex-valued function $f$ on $F^n$, we define the \emph{Fourier transform} 
$\hat f$ on $F^n$ by
$$ 
\hat f(\xi) := \int_{F^n}  f(x)\ e(-x \cdot \xi)\ dx
$$
where $x \cdot \xi$ is the inner product
$$ 
(x_1, \ldots, x_n) \cdot (\xi_1, \ldots \xi_n) := x_1 \xi_1 + \ldots + x_n \xi_n.
$$
The inverse of this transform is given by
$$
g \spcheck(x) := \int_{F^n_*} g(\xi)\ e(x \cdot \xi)\  d\xi.$$
We have the Parseval identity
$$ 
\int_{F^n} f_1(x) \overline{f_2(x)}\ dx = \int_{F^n_*}
\hat f_1(\xi) \overline{\hat f_2(\xi)}\ d\xi,
$$
the Plancherel identity
$$
\|f\|_{L^2(F^n, dx)} = \|\hat f\|_{L^2(F^n_*, d\xi)},
$$
and the intertwining of convolution and multiplication
$$ \hat f_1 \hat f_2 = \widehat{f_1 * f_2}; \quad
\widehat{f_1 f_2} = \hat f_1 * \hat f_2.$$
Note that $f_1 * f_2$ is convolution using counting measure $dx$, while 
$\hat f_1 * \hat f_2$ is convolution using normalized counting measure $d\xi$.

\section{The restriction problem}\label{rest-sec}

We begin by setting up some notation for the restriction problem.

Let $S \subset F^n_*$ be a non-empty set of frequencies, which in this 
paper will always be an algebraic variety in $F^n_*$.  We endow $S$ with 
normalized ``surface measure''
$$ 
\int_S g(\omega)\ d\sigma(\omega) := 
\frac 1{|S|} \sum_{\xi \in S}\ g(\xi),
$$
thus $S$ has total mass 1 and by identifying measures with functions by dividing 
by the ambient measure we have
$$
(g d\sigma)\spcheck(x) = \frac 1{|S|}\sum_{\xi\in S}\ g(\xi)\ e(x\cdot \xi).
$$
For any exponents $1 \leq p, q \leq \infty$, we let $R^*(p \to q) = R^*_S(p \to q)$ 
be the best constant such that the estimate
\be{rest-def}
\| (g d\sigma)\spcheck \|_{L^q(F^n, dx)}\, \leq\, R^*(p \to q) \ 
\| g \|_{L^p(S, d\sigma)}.
\end{equation}
holds for all functions $g$ on $S$.  Since $F$ is finite, the quantity 
$R^*(p \to q)$ will be a finite positive real.  By duality, one can also 
define $R^*(p \to q)$ to be the best constant such that we have the 
restriction estimate
\be{rest-dual}
\| \hat{f} \|_{L^{p'}(S, d\sigma)} \,\leq \,R^*(p \to q) \ \| f \|_{L^{q'}(F^n, dx)}
\end{equation}
for all functions $f$ on $F^n$.

For instance, it is immediate from the Plancherel's identity and H\"older's 
inequality that
\be{triv}
R^*(p \to \infty) = 1 \ \ \hbox{ for }\ 1 \leq p \leq \infty
\end{equation}
and
\be{triv-2}
R^*(p \to 2) = \left(\frac{{}\ |F|^n}{|S|}\right)^{ 1/2} \ \ \hbox{ for }\  2 \leq p \leq \infty .
\end{equation}

For any algebraic variety $S$, i.e. a zero set in $F^n_*$ of a polynomial in $\Z[\xi_1,\dots,\xi_n]$, 
the \emph{restriction \- problem} for $S$ 
asks to determine the set of exponents $p$ and $q$ such that 
$R^*(p \to q) \leq C_{p,q}$, where $C_{p,q}$ does not depend on the 
underlying field $F$.  One may think of the restriction problem as a 
quantitative way to understand the average size of exponential sums 
of the form
$$ 
\sum_{\xi \in S}\ g(\xi)\ e(x \cdot \xi)
$$
where $g$ is an arbitrary  sequence in $l^p(S)$.  To 
put it loosely, restriction problems are an attempt to understand 
exponential sums in which the coefficients $g(\xi)$ have no easily 
exploitable structure other than magnitude bounds. 

When the underlying field $F$ is replaced by $\R$, and $S$ is the unit sphere  
in $\R^n_*$, this problem was posed by E. M. Stein in 1967 
\cite{feff:thesis}, who observed (in our notation) that $R^*(p \to q)$ 
can be bounded for certain $q < \infty$ if $S$ has some non-trivial 
curvature. 

To see how the Euclidean and the finite field restriction estimates 
are related consider the following heuristic argument 
(by standard arguments this can be made rigorous):
Let  $M\subset R^n$ be a compact $d$-dimensional smooth submanifold 
such that the translations of  $M$ along normal vectors $n\neq 0$ 
at some point of $M$ do not intersect for $0<|n|<1/R$ and 
suppose $\supp(f)\subset Q_R$, a cube of sidelength 
$\approx R$ around the origin. By the uncertainty principle
$\hat f$ is roughly constant on a ball of radius $1/R$. 
Hence, if $M_R$ is a $1/R$-neighborhood of $M$ we have with 
$$
A=\{m\in Z^n|\ \dist(\frac mR , M) \leq \frac 1R\}
       \qquad \text{(note that}\  |A|\approx R^d\text{)},
$$
\begin{align*}
\int_{M_R} \ |\hat f|^{p'}\  dx
&\approx \sum_A 
\ \int_{B(\frac mR, \frac 1R)}\  |\hat f(x)|^{p'} \ dx \\
&\approx  \sum_A \ \Vol_n(B_{\frac 1R}(\frac mR))\ \  
          |\hat f(\frac mR)|^{p'}\\
&= \frac {c_n}{R^n} \ \sum_A \  \  |\hat f(\frac mR)|^{p'}.
\end{align*}
Let $f_R(x)=f(Rx)$, thus $f_R$ has support in the unit cube and 
$\hat f(\frac mR)=R^n\ \hat f_R(m)$. Hence 
$$
\int_{M_R}\ \ |\hat f|^{p'}\  dx \ \approx\  \frac {R^{np'}}{R^n}\ \sum_A \  \  |\hat f_R(m)|^{p'}.
$$
On the other hand we may assume $M_R \subset \bigcup_{|nR|<1}  (n+M)$. Hence
$$
\int_{M_R} |\hat f|^{p'}\ dx \lesssim  
		\int_{n\in R^{n-d}, |n|<\frac 1R}\ \int_{n+M} \ |\hat f|^{p'}.
$$
By translation invariance and assuming a $(q',p')$-restriction inequality for $M$:
\begin{align*}
\int_{M_R} |\hat f|^{p'} \ dx    & 
	\lesssim  \ \frac 1{R^{n-d}}\ \sup_n\ \int_{n+M} \ |\hat f|^{p'}\  dx \\
&\lesssim \ \frac 1{R^{n-d}}\ \|f\|_{L^{q'}(Q_R)}^{p'} \ = 
\frac { R^{ \frac {np'}{q'}}} {R^{n-d}}\ \left( \int_{x\in Q_1}\ |f_R(x)|^{q'}\ dx\right)^{p'/q'}.
\end{align*}
With $g=f_R$ considered as a periodic function with Fourier coefficients $\hat g(m)$ we find:
$$
 \frac {R^{np'}}{R^n}\ \sum_A \  \  |\hat g(m)|^{p'} \lesssim \ \frac { R^{ \frac {np'}{q'}}} {R^{n-d}}\    \left(  \int_{x\in Q_1}\ |g(x)|^{q'} \ dx\right)^{p'/q'}
$$
i.e.
$$
\Big(\sum_A \  \  |\hat g(m)|^{p'}\Big)^{\frac 1{p'}} \lesssim
\ R^{\frac d{p'}} \ R^{-\frac n{q}}\  \|g\|_{L^{q'}(T^n)}.
$$
By duality we get 
$$
\Big(\int_{T^n}\ \big| \sum_{m\in A} \  a_m  \ e^{im\cdot x} \big|^{q}\ \Big)^{\frac 1{q}}
\lesssim
 \ R^{\frac d{p'}} \ R^{-\frac n{q}}\ 
\Big( \sum_{m\in A} \  |a_m|^{p}\Big)^{\frac 1{p}}. 
$$
By eventually translating $A$ we may assume that $A\subset [0,R)^n\cap Z^n$. 
If we discretize the integral by a Riemann sum over cubes of sidelength $1/R$ 
according to a result of Marcinkiewicz and Zygmund (see \cite[Vol II p. 28]{zyg}) 
the $L^q$-norm on the right hand side is equivalent to the discretized 
sum $(1<q<\infty)$:
$$
\Big( \frac 1{R^n} \sum_{0\leq x_i<R, x_i\in \Z} \big| \sum_{m\in A} \ 
                   a_m  \ e^{i m\cdot \frac{x}{R}}\ \big|^{q}\ \Big)^{\frac 1{q}}
\lesssim
 \ R^{-\frac n{q}}\ R^d\ R^{-\frac d{p}} \ 
\Big( \sum_{m\in A} \  |a_m|^{p}\Big)^{\frac 1{p}}. 
$$		   
Therefore, if $p\approx R$ is prime and $F=\Z/p\Z$ by identifying $A$ 
with a subset $S$ in $F_*^n$ we get, since $|S|\approx R^d$:
$$
\Big(  \sum_{x\in F^n} \big| \frac 1{|S|}\ \sum_{m\in S} \ 
                   a_m  \ e(m\cdot x)\ \big|^{q}\ \Big)^{\frac 1{q}}
\lesssim \ 
 \ 
\Big( \frac 1{|S|}\  \sum_{m\in S} \  |a_m|^{p}\Big)^{\frac 1{p}},
$$
i.e. 
$$
R_S^*(p\to q) \lesssim 1.
$$
Note that $S\subset F^n_*$ depends on the prime $p$, i.e.
here $S$ is not given by a fixed algebraic equation over $\Z$.

A large amount of work, especially in the three-dimensional case, 
has since been done on the restriction problems 
(see e.g. \cite{stein:large}, \cite{borg:stein}, \cite{wolff:cone}, 
\cite{tv:cone1} and the references therein); notably for the case of a 
circle in $\R^2$ the problem was solved by C. Fefferman and E. M. Stein 
in \cite{feff:thesis}. The restriction problem is also connected to 
problems in dispersive and wave equations, and to 
certain problems in number theory; we refer the reader to \cite{Bo}, 
\cite{wolff:survey}, \cite{mock:habil} for further 
discussion of these connections.

In this paper we shall mostly restrict ourselves to the case when $S$ is a 
hypersurface, however the lower codimension case is also extremely 
interesting (see \cite{christ:thesis}, \cite{mock:habil}).  
In fact, we shall mostly concern ourselves with paraboloids and cones in $F^n_*$; 
we have avoided the sphere as the Fourier transform of surface 
measure is not as easy to compute in the finite field case.  Other 
quadratic surfaces are also of interest, see e.g. 
\cite{strichartz:restrictionquadratic}.

\section{General restriction theory}\label{basic-sec}

In this section we collect some general (and rather easy) facts about the 
quantity $R^*(p \to q)$.  We shall assume that $|S| \sim |F|^d$ for 
some $0 < d < n$ (this is the finite field analogue of $S$ having 
dimension $d$).

From H\"older's inequality and the embedding $l^{p_1} \subset l^{p_2}$ for 
$p_1 \geq p_2$ we see that 
\be{h1}
R^*(p_1 \to q) \leq R^*(p_2 \to q) \hbox{ when } p_1 \geq p_2
\end{equation}
and
\be{h2}
R^*(p \to q_1) \leq R^*(p \to q_2) \leq |F|^{n(\frac{1}{q_2} - \frac{1}{q_1})}\ 
R^*(p \to q_1) \hbox{ when } q_1 \geq q_2.
\end{equation}
From these inequalities and \eqref{triv-2} we see that
\be{junk}
R^*(p \to q) \geq |F|^{-n (\frac{1}{2} - \frac{1}{q})_+}\ R^*(\infty \to 2)
= |F|^{-n(\frac{1}{2}-\frac{1}{q})_+}\ \Big(\frac{{}\ |F|^n}{|S|}\Big)^{1/2}.
\end{equation}
Here we use $a_+$ as shorthand for $\max(a,0)$.
Since $|S| \sim |F|^d$, we thus see that $R^*(p \to q)$ can only be 
bounded by $O(1)$ if
\be{q-1}
q \geq \frac{2n}{d}.
\end{equation}

From \eqref{junk} we have in particular that
$$ 
R^*(2 \to q) \gtrsim |F|^{n/q}\ |S|^{-1/2}
$$
whenever $2 \leq q \leq \infty$.  If the $\gtrsim$ can be replaced 
by $\sim$, then $R^*(2\to q)$ resembles the $\Lambda(q)$-constant
of the set $S$ 
(see e.g. \cite{borg:lambda-p}); roughly speaking, this means that the 
exponentials $\{ e(x \cdot \xi): \xi \in S\}$ form an essentially 
orthogonal set in $L^q$.

Also, if we test \eqref{rest-def} with $g$ equal to a Dirac delta 
$g = \delta_\eta$ for some $\eta \in S$, we obtain
$$ 
|S|^{-1}\ |F|^{n/q} \leq R^*(p \to q)\ |S|^{-1/p}
$$
which implies that $R^*(p \to q)$ can only be bounded by $O(1)$ if
\be{q-2}
q \geq \frac{np'}{d}.
\end{equation}

As we shall see, in some cases the necessary conditions \eqref{q-1}, 
\eqref{q-2} for the boundedness of $R^*(p \to q)$ are also sufficient.  
However, if $S$ is somehow ``flat'', then one can improve the above 
necessary conditions.

For instance, suppose $S$ contains an affine subspace $V \subset F^n_*$ of 
dimension $k$.  Then if we test \eqref{rest-def} with $g$ equal to the 
characteristic function $\chi_V$, we obtain
$$ |F|^k\ |S|^{-1}\ |F|^{(n-k)/q} \leq R^*(p \to q)\ \, |F|^{k/p}\ |S|^{-1/p}$$
which yields the condition
\be{q-improv}
q \geq p' \frac{n-k}{d-k}.
\end{equation}

The problem of estimating $R^*(p \to q)$, in particular 
estimates for $R^*(\infty\to q)$, is related to the Hardy-Littlewood 
majorant problem, see \cite{HL,Li}.  This problem asks 
whether for any $2\leq p < \infty$ and any functions $f, g$ 
with $|\widehat f(\xi)| \leq \widehat g(\xi)$ for all 
$\xi \in F^n_*$, one has $\|f\|_p \leq C_p(F)\, \|g\|_p$ with $C_p(F)$
independent of $F$.
It was observed by Hardy and Littlewood in \cite{HL} 
that this is true with $C_p(F)=1$ for even integer $p$, 
but examples of Boas, Bachelis, and Fournier show that $C_p(F)$ 
grows unboundedly with $F$ for all other values of $p$ 
(see e.g. \cite{fournier}).  
A quantitative lower bound for the grows of  $C_p(F)$ is given in \cite{mock:habil},
e.g. for $2<p<4$ one has $C(F)\geq  |F|^{\frac {c_p}{\log\log |F|}}$ 
with $c_p>0$ independent of $F$.
However, this failure is only $O(|F|^\eps)$ for any $\eps > 0$, 
so it is still possible that a weaker version 
$\|f\|_p \lessapprox \|g\|_p$ holds.  
If this were true then we would have in 
particular that (for all $S$)
$$ 
\R^*(\infty \to q) \approx \| (d\sigma) \spcheck \|_q.
$$ For further discussion see \cite{mock:habil}.

\section{Even exponents}\label{even-sec}

When $q$ is an even integer one can compute $R^*(p \to q)$ quite 
easily just from Plancherel's theorem.  For instance:

\begin{lemma}\label{k-planch}  Let $q = 2k$ be a positive even integer, 
and suppose that for any $\eta \in F^n_*$ the number of solutions to the 
equation
\be{sum}
\eta = \xi_1 + \ldots + \xi_k; \quad \xi_1, \ldots, \xi_k \in S
\end{equation}
is bounded by $A$.  Then we have
$$R^*(2 \to 2k) \leq A^{1/2k}\ |F|^{n/2k}\ |S|^{-1/2}.$$
\end{lemma}

\begin{proof}
We have to show that
$$ 
\| \widehat{f d\sigma} \|_{L^{2k}(F^n, dx)} \leq 
A^{1/2k}\ |F|^{n/2k}\ |S|^{-1/2}\ \| f \|_{L^2(S, d\sigma)}.
$$
Raising this to the $k^{th}$ power and applying Plancherel, this becomes
$$ 
\| f d\sigma * \ldots * f d\sigma \|_{L^{2}(F^n_*, d\xi)} \leq 
A^{1/2}\ |F|^{n/2}\ |S|^{-k/2}\ \| f \|_{L^2(S, d\sigma)}^k
$$
where the convolution contains $k$ copies of $f d\sigma$.  

It is now clear that we can assume $f$ is positive.
By hypothesis we have
$$ 
\| d\sigma * \ldots * d\sigma \|_{L^\infty(F^n_*, d\xi)} \leq A\ |F|^n\ |S|^{-k}
$$
while from Fubini's theorem we have
$$ 
\| f^2 d\sigma * \ldots * f^2 d\sigma \|_{L^1(F^n_*, d\xi)} = 
\| f \|_{L^2(S, d\sigma)}^{2k}.
$$
The claim then follows by taking the geometric mean of these estimates 
and applying Cauchy-Schwartz.
\end{proof}

The above lemma has seen many applications, e.g. Sidon 
employed it in his study of lacunary Fourier series.
A variant in Euclidean spaces using Hausdorff-Young's inequality 
instead of Plancherel's identity is the crucial tool in settling 
the restriction problem for circles and more general curves in $\R^2$ 
(see e.g. \cite{feff:thesis}, \cite{sj}).  The continuous 
variants as well as related estimates in the theory of bilinear  
$X^{s,b}$-estimates are also useful in non-linear PDE's 
(see e.g. \cite{borg:xsb}, \cite{klainerman:nulllocal}, \cite{tao:multi}, etc.).

As an application we consider the polynomial curve
$$ 
S := \gamma(F),
$$
where
$$ 
\gamma(t) := (t,t^2, \ldots, t^n).
$$
The $n$ coordinates of $\gamma(t_1) + \ldots + \gamma(t_n)$ are just the 
first $n$ symmetric functions of $t_1, \ldots, t_n$.  If $\chr (F) > n$ Newton's
formulae show that the $t_1, \ldots, t_n$ are determined up to permutations.  
Thus we can apply the above Lemma with $k := n$ and $A := n!$ to conclude that 
$R^*(2 \to 2n)$ is bounded.  From \eqref{h1}, \eqref{h2} we thus 
see that $R^*(p \to q)$ is bounded whenever
$$ 
q \geq 2n, np'.
$$
This exactly matches the necessary conditions \eqref{q-1}, \eqref{q-2} 
(with $d=1$).  Note that if $F$ has characteristic $n$ or less, then 
$\gamma$ becomes degenerate and $S$ is contained inside a proper 
subspace of $F^n_*$.

For the rest of the section we consider the paraboloid
\be{paraboloid}
S := \{ (\xi, \xi \cdot \xi): \xi \in F_*^{n-1} \}.
\end{equation}

We apply Lemma \ref{k-planch} with $k := 2$.  To find $A$, we need to 
bound solutions to the problem
$$ 
\xi_1 + \xi_2 = \eta', \quad \xi_1 \cdot \xi_1 + \xi_2 \cdot \xi_2 = \tau
$$
for any $\eta' \in F^{n-1}$ and $\tau \in F$.  
If one fixes $n-2$ of the coordinates of $\xi_1$, then the remaining 
coordinate obeys a quadratic equation (after substituting 
$\xi_2 = \eta' - \xi_1$) and thus has at most two solutions.  Thus we may 
take $A = 2 |F|^{n-2}$.  From Lemma \ref{k-planch} we thus see that 
$R^*(2 \to 4)$ is bounded.
(The same argument gives $R^*(2 \to 2k)$ for any even integer $2k$, but 
this is in any event implied from the $2k=4$ case by \eqref{h2}).

From \eqref{q-1}, \eqref{q-2} with $d := n-1$ we obtain the necessary 
conditions
$$ 
q \geq \frac{2n}{n-1}, \frac{np'}{n-1}
$$
for the boundedness of $R^*(p \to q)$.  If $n=2$ then the above discussion 
(combined with \eqref{h1}, \eqref{h2}) shows that these necessary conditions 
are in fact sufficient; this was first observed by Zygmund.  

For $n > 2$ it is possible for the paraboloid to contain lines 
(since $\xi \cdot \xi$ can vanish even when $\xi \neq 0$) or even 
higher-dimensional subspaces, at which point one can improve the above 
necessary conditions.  To avoid these complications let us restrict 
ourselves to the case $n=3$, with $-1$ not a square number in $F$ 
(so that $|F|$ is a power of a prime $p = 3 (\mod 4)$).  In this case the 
parabola does not 
contain any lines, and the necessary conditions become
$$ 
q \geq 3, \frac{3p'}{2}.
$$
It seems reasonable to conjecture that the above conditions are still 
sufficient.  In other words, we conjecture that $R^*(2 \to 3)$ is 
bounded - which is substantially stronger than the bound of 
$R^*(2 \to 4)$ obtained earlier.  We remark that the corresponding 
conjecture for the paraboloid in $\R^3$ is a challenging open question; see 
e.g. \cite{tv:cone2} for some recent progress.

As partial evidence of this conjecture we present

\begin{proposition}\label{12-7}  If $n=3$, $-1$ is not a square number, 
and $S$ is the paraboloid \eqref{paraboloid}, then we have
$$
R^*(8/5 \to 4) \lessapprox 1.
$$
\end{proposition}

Note that this is on the line $q = 3p'/2$.  This only improves the 
``$p$'' index of the $R^*(2 \to 4)$ result, but we will give an 
improvement to the ``$q$'' index later in Theorem \ref{BR}.  The 
logarithmic factor can probably be removed, but we will not do so here.

Before we prove Proposition \ref{12-7} we need the following 
well-known incidence geometry lemma (see e.g. \cite{bollobas}).

\begin{proposition}\label{lines}  Let $P$ be a collection of points 
in $F^2$, and let $L$ be a collection of lines in $F^2$.  Then
$$
|\{ (p, l) \in P \times L: p \in l \}| \leq 
\min( |P|^{1/2}\, |L| + |P|, |P|\, |L|^{1/2} + |L|)
$$
\end{proposition}

\begin{proof} By the duality of points and lines in $F^2$ (or more precisely, 
projective space $PF^2$) it suffices to verify the first bound.  
Denote the set on the left-hand side by $I$.  By \eqref{cauchy} 
we have 
$$ 
|\{ (p,l,l') \in P \times L \times L: p \in l; p \in l' \}|
\geq \frac{|I|^2}{|P|}.
$$
Fix $l, l'$.  If $l \neq l'$ then there is at most one point $p$ 
which contributes to the left-hand side.  Thus
$$ 
|I| + |L|^2 - |L| \geq \frac{|I|^2}{|P|}
$$
which we re-arrange as
$$ 
(\, |I| - |P|/2\, )^2 \leq |P|\, |L|^2 - |P|\, |L| + |P|^2/4
\leq (\ |P|^{1/2}\, |L| + |P|/2 \ )^2.
$$
The claim follows.
\end{proof}

We can now prove Proposition \ref{12-7}.

\begin{proof}  This will be  analogous to the ``12/7'' estimate in 
\cite{vargas:restrict}, \cite{tvv:bilinear}.  
We want to show that
$$ 
\| (f d\sigma)\spcheck \|_{L^4(F^n_*, d\xi)} \lessapprox 
\| f \|_{L^{8/5}(S, d\sigma)}.
$$
Since we are allowing ourselves to lose logarithmic powers of $|F|$ 
we may restrict ourselves to characteristic functions $f = \chi_E$ 
by the usual dyadic pigeonholing argument, where $E$ is any non-empty 
subset of $S$.

Fix $E$.  By arguing as in Lemma \ref{k-planch}, it suffices to show that
$$ 
\| \chi_E d\sigma * \chi_E d\sigma \|_{L^2(F^n_*, d\xi)} \lesssim 
\| \chi_E \|_{L^{8/5}(S, d\sigma)}^2.
$$
Expanding this all out, this becomes
$$ 
\sum_{\omega_1, \omega_2, \omega_3, \omega_4 \in E: \omega_1 + \omega_2 = 
\omega_3 + \omega_4} 1 \lesssim |E|^{5/2}.
$$
Since $E \subset S$, it will suffice to show that
$$ 
\sum_{\omega_2, \omega_3 \in E: \omega_3 - \omega_2 \in S - \omega_1} 1
\lesssim |E|^{3/2}
$$
for all $\omega_1 \in E$.
Since $S$ is invariant under Galilean transformations
$$ 
g_a:F_*^{n-1}\times F_*\ni (\xi, \tau) \mapsto 
(\xi + a, \tau + 2 \xi \cdot a +  a\cdot a)
$$
 the above sum transforms under $g_{-\eta_1}$, where 
$(\eta_1, \eta_1\cdot \eta_1)=\omega_1$, into  
$$ 
\sum_{\omega_2, \omega_3 \in E': \omega_3 - \omega_2 \in S} 1
$$
with $E'=g_{-\eta_1}(E)$. 
We may take $\omega_2, \omega_3 \neq 0$ since the contribution of 
$\omega_2 = 0$ or $\omega_3 = 0$ is $O(|E|)$.  If we write 
$E' \backslash \{0\} = \{ (\xi, \xi \cdot \xi): \xi \in P \}$ for 
some $P \subset F_*^2$, and change coordinates such that 
$\omega_2 = (\xi, \xi \cdot \xi)$ and $\omega_3 = (\eta,\eta \cdot \eta)$, 
the condition $\omega_3 - \omega_2 \in S$ then becomes 
$\xi \cdot \eta = \xi \cdot \xi$, and we reduce to showing
$$
|\{\, (\xi, \eta) \in P \times P: \xi \cdot \eta = \xi \cdot \xi \,\}| 
\lesssim |P|^{3/2}.
$$
For any $\xi \in P$, let $l(\xi)$ denote the line
$$ 
l(\xi) := \{ \eta \in F_*^2: \xi \cdot \eta = \xi \cdot \xi \}
$$
and let $L$ denote the set of all such lines $l(\xi)$.  Since $-1$ 
is not a square we see that these lines are all distinct,
so $|L| = |P|$.   Our task is then to show the incidence bound
$$
|\{ (p, l) \in P \times L: p \in l \}| \lesssim |P|^{3/2}.
$$
But this follows from Proposition \ref{lines}.
\end{proof}

\section{Tomas-Stein type arguments}\label{tomas-sec}

Let $S$ be any subset of $F^n_*$.  Suppose we wish to estimate the 
quantity $R^*(2 \to q)$ for some $q$.  If we square \eqref{rest-dual} 
and apply Plancherel, we see that
\be{tts}
|\langle f, f * (d\sigma)\spcheck \rangle| \leq R^*(2 \to q)^2\
\| f \|_{L^{q'}(F^n, dx)}\ \| f\|_{ L^{q'}(F^n, dx)}
\end{equation}
and furthermore that $R^*(2 \to q)$ is the best constant with this property.

It is thus of interest to study the object $(d\sigma)\spcheck$.  At the 
origin we clearly have $(d\sigma)\spcheck(0) = 1$, however away from the 
origin we often have some decay.  To this end we introduce the 
\emph{Bochner-Riesz kernel}
$$ K(x) := (d\sigma)\spcheck(x) - \delta_0(x)$$
which is just $(d\sigma)\spcheck$ with the origin removed.  This 
kernel is the finite field analogue of (a dyadically localized 
portion of) the kernel for Bochner-Riesz means in $\R^n$, see e.g. 
\cite{feff:thesis}, \cite{carl:disc} for further discussion.

Suppose that we have the decay estimate
\be{decay-ass}
\| K \|_{L^\infty(F^n, dx)} \lesssim |F|^{-\tilde d/2}
\end{equation}
for some exponent $0 < \tilde d < n$; this quantity is sometimes 
referred to as the \emph{Fourier dimension} of $S$.

If $|S| \sim |F|^d$ for some $0 < d < n$, then from Plancherel 
we have
$$ \| (d\sigma)\spcheck \|_{L^2(F^n, dx)} = \| d\sigma \|_{L^2(F^n_*, d\xi)}
= |F|^{(n-d)/2},
$$
which implies that $\tilde d \leq d$.  In some cases (e.g. the 
parabola) the two notions of dimension are equal, but when $S$ 
has some arithmetic closure properties the Fourier dimension is 
often smaller than the ordinary dimension.

If we could somehow ignore the $x = 0$ component of $(d\sigma)\spcheck(x)$, 
then from \eqref{decay-ass}, then we would have
$$ 
|\langle f, f * (d\sigma)\spcheck \rangle| \lesssim |F|^{-\tilde d/2}\
\| f \|_{L^1(F^n, dx)}\ \| f\|_{ L^1(F^n, dx)}
$$
which would then imply $R^*(2 \to \infty) \lesssim |F|^{-\tilde d/4}$.  
Of course this is absurd by \eqref{triv}, however for interpolation 
purposes this is morally true\footnote{This idea of passing to a 
bilinear formulation and eliminating the diagonal case is quite 
useful in the Euclidean theory, see e.g. \cite{tvv:bilinear}, 
\cite{tv:cone1}, \cite{tv:cone2}, \cite{wolff:cone}.}.  A precise 
statement is the following (which is closely related to the 
formulation of the Tomas-Stein argument given in Bourgain \cite{borg:kakeya}):

\begin{lemma}[Tomas-Stein argument]\label{ts}  Let $p, q \geq 2$, 
and suppose that \eqref{decay-ass} holds.  Then for any $0 < \theta < 1$ 
we have
$$ 
R^*(p \to q/\theta) \lesssim 1 + R(p \to q)^\theta\ |F|^{-\tilde d(1-\theta)/4}.
$$
\end{lemma}

\begin{proof}
Fix $f$ to be in the unit ball of $L^{(q/\theta)'}(F^n, dx)$.  It 
suffices to show that
\be{ts-targ}
\| \hat f \|_{L^{p'}(S, d\sigma)} \lesssim 
1 + R^*(p \to q)^\theta\ |F|^{-\tilde d(1-\theta)/4}.
\end{equation}
Let $\lambda > 0$ be a parameter to be chosen later.  We can divide 
into two cases: either $|f(x)| \leq \lambda$ for all $x \in F^n$, 
or else $|f(x)| \geq \lambda$ for all $x$ in the support of $f$.

In the first case we have
$$ 
\| f \|_{L^{q'}(F^n, dx)} \lesssim \lambda^{(q' - (q/\theta)')/q'}
= \lambda^{-(1-\theta)/(q-\theta)}
$$
and so by \eqref{rest-dual} we have
$$
\| \hat f \|_{L^{p'}(S, d\sigma)} \lesssim R^*(p \to q) \
\lambda^{-(1-\theta)/(q-\theta)}.
$$

Now suppose we are in the second case.  We can estimate the 
$L^{p'}(S, d\sigma)$ norm by the $L^2(S, d\sigma)$ norm and use 
Plancherel to obtain
$$
\| \hat f \|_{L^{p'}(S, d\sigma)}^2 \leq 
|\langle f, f * (d\sigma)\spcheck \rangle|.
$$
We split $(d\sigma)\spcheck = \delta_0 + K$.  Using \eqref{decay-ass} we obtain
$$
\| \hat f \|_{L^{p'}(S, d\sigma)}^2 \leq 
\| f\|_{L^2(F^n, dx)}^2 + |F|^{-\tilde d/2}\ \| f \|_{L^1(F^n, dx)}^2.$$
Since $dx$ is discrete we have
$$\| f\|_{L^2(F^n, dx)} \leq \| f \|_{L^{(q/\theta)'}(F^n, dx)} \leq 1$$
while by the assumption $|f(x)| \geq \lambda$ on the support of $f$ we have
$$ \| f \|_{L^1(F^n, dx)} \lesssim \lambda^{1-(q/\theta)'}
= \lambda^{-\theta/(q-\theta)}.$$
Thus we have
$$ 
\| \hat f \|_{L^{p'}(S, d\sigma)} \lesssim 
1 + |F|^{-\tilde d/2}\, \lambda^{-\theta/(q-\theta) }.
$$
To balance our two cases we define $\lambda > 0$ by
$$ 
R^*(p \to q) \lambda^{-(1-\theta)/(q-\theta)} =
|F|^{-\tilde d/2}\, \lambda^{-\theta/(q-\theta) }.
$$
A little algebra then shows that \eqref{ts-targ} holds in both 
cases\footnote{Alternatively, one can proceed by considering 
convolution with $\delta_0 + |F|^z K$ where $z$ is a complex 
interpolation parameter; this is in the spirit of Stein's argument 
\cite{stein:large}. 
We omit the details.}.
\end{proof}

In particular, if we insert the bound $R^*(2 \to 2) \lesssim |F|^{(n-d)/2}$ 
from \eqref{triv-2} we obtain that $R^*(2 \to q)$ is bounded whenever
$$ 
q \geq 2 + \frac{4 (n-d)}{\tilde d}.
$$

In the Euclidean setting this bound on $q$ is often quite sharp (if $p$ is 
constrained to equal 2).  This is usually due to the presence of ``Knapp'' 
examples for $S$, see e.g. \cite{tomas:restrict}.  The analogous type of 
example in the Euclidean case would be an arithmetic progression, or perhaps 
a subspace as in \eqref{q-improv}, but these are not always available.

To illustrate this let us again take the paraboloid \eqref{paraboloid}.   
If $x = (x_1, \ldots, x_n)$ with $x_n \neq 0$, then
\be{gauss}
\begin{split}
(d\sigma)\spcheck(x) &= |F|^{1-n} \sum_{\xi \in F_*^{n-1}} 
e(x \cdot \xi + x_n\, \xi \cdot \xi)\\
&= |F|^{1-n} \prod_{i=1}^n \sum_{\xi_i \in F_*}
e(x_i \xi_i + x_n \xi_i \xi_i)\\
&= |F|^{1-n} \prod_{i=1}^n 
e(x_i x_i / 4 x_n) \sum_{\xi_i \in F_*}
e(x_n (\xi_i + x_i/2x_n)^2)\\
&= |F|^{1-n} S(x_n)^{n-1}
e(\underline{x} \cdot \underline{x} / 4 x_n)
\end{split}
\end{equation}
where $\underline{x} := (x_1, \ldots, x_{n-1})$ and $S(x)$ is the Gauss sum
$$ 
S(x) := \sum_{\xi \in F_*} e(x \xi^2).
$$
If $x_n = 0$, then the above sum collapses instead to
$$ 
(d\sigma)\spcheck(x) = \delta_0(\underline{x}).
$$
As is well known, the magnitude of the Gauss sum is $|F|^{1/2}$.  Indeed
we have
\be{gauss-size}
\begin{split}
|S(x)|^2 &= \sum_{\xi, \eta \in F_*} e(x \xi^2)\ e(- x\eta^2) \\
&= \sum_{\xi, \eta \in F_*} e(x (\xi+\eta) (\xi-\eta)) \\
&= |F| + \sum_{\xi, \eta \in F_*: \xi \neq \eta} e(x (\xi+\eta) (\xi-\eta)) \\
&= |F| + \sum_{h \in F_*: h \neq 0} \sum_{\eta \in F} e(x (2\eta+h) h) \\
&= |F| + \sum_{h \in F_*: h \neq 0} 0.
\end{split}
\end{equation}
Thus \eqref{decay-ass} holds with $\tilde d := d = n-1$.  In particular 
we obtain that $R^*(2 \to q)$ is bounded whenever
\be{ts-exp}
q \geq 2 + \frac{4}{n-1}.
\end{equation}
This does not compare well with the results of the previous section, 
although it does reproduce the result $R^*(2 \to 4)$ proved earlier 
in the case when $n=3$ and $-1$ is not a square.  The condition 
\eqref{q-2} (or \eqref{q-1}) suggests that one should in fact obtain 
boundedness of $R^*(2 \to q)$ whenever
$$ 
q \geq 2 + \frac{2}{n-1},
$$
which is indeed the case in two dimensions.  However, in higher dimensions 
the paraboloid can contain subspaces in which case one can narrow the gap.  
For instance, when $n=3$ and $-1 = i^2$ for some $i \in F$, then $S$ 
contains lines such as $\{ (t, it, 0): t \in F \}$ and so by 
\eqref{q-improv} we obtain a necessary condition of $q \geq 4$, which 
does indeed match the Tomas-Stein exponent \eqref{ts-exp} in this case.  
However, this does not preclude that the bound of $R^*(2 \to 4)$ can be 
improved in other directions, indeed from \eqref{q-1}, \eqref{q-improv} 
one might conjecture that $R^*(3 \to 3)$ holds.  (An analogous conjecture 
exists for the paraboloid in $\R^3$).

Similar examples can be made in higher odd dimensions (with $k := (n-1)/2$), 
but there appears to be a gap between the Tomas-Stein exponent and the best 
known counterexamples when $n$ is even or when $-1$ is  not a square.  We 
shall not pursue these matters further, however we shall close this section 
with an improvement of the ``$q$'' index of the Tomas-Stein restriction 
theorem $R^*(2 \to 4)$ in the case when $n=3$ and $-1$ is not a square.

\begin{theorem}\label{BR}
When $n=3$, $S$ is the paraboloid \eqref{paraboloid} and $-1$ is not a 
square, then $R^*(2 \to 18/5 + \eps) \leq C_\eps$ for all $\eps > 0$.
\end{theorem}

\begin{proof}
We shall exploit the restriction theorem in Proposition \ref{12-7}
as well as Fourier identities for the paraboloid which in the Euclidean setting 
have been exploited in \cite{carbery:parabola}.  We also use the 
close relationship between the restriction and Bochner-Riesz problems, 
first observed in \cite{feff:thesis}.

By Lemma \ref{ts} with $p=2$, $q = 16/5$, $\theta = 8/9-$ it suffices to 
show that
$$ 
R^*(2 \to 16/5) \lessapprox |F|^{1/16}.
$$
(Note from \eqref{junk} that this is sharp up to logarithms).
By \eqref{tts} it suffices to show that
$$ 
|\langle f, g * (d\sigma)\spcheck \rangle| 
\lessapprox |F|^{1/8}\ \| f \|_{L^{16/11}(F^3, dx)}\ \| g\|_{L^{16/11}(F^3, dx)};
$$
for all $f, g$ on $F^3$. It suffices to show that
$$ 
|\langle f, g * (d\sigma)\spcheck \rangle| \lessapprox |F|^{1/8}\
\| f \|_{L^{4/3}(F^3, dx)}\ \| g\|_{ L^{8/5}(F^3, dx)},
$$
since bilinear interpolation of this inequality with 
the one where $f$ and $g$ are exchanged gives 
the above estimate on $L^{16/11}$. 
By duality above inequality becomes
$$ 
\| g * (d\sigma)\spcheck \|_{L^4(F^3, dx)} \lessapprox |F|^{1/8}\ \| g \|_{L^{8/5}(F^3, dx)}.
$$
The $\delta_0$ portion of $(d\sigma)\spcheck$ is definitely acceptable by 
Young's inequality, so it suffices to show
$$
\| g * K \|_{L^4(F^3, dx)} \lessapprox |F|^{1/8}\ \| g \|_{L^{8/5}(F^3, dx)}$$
where $K$ is the Bochner-Riesz kernel introduced earlier.
We split $x := (x_1, x_2, x_3)$, and for each $y_3 \in F$ define the 
function $g_{y_3}$ to be the restriction of $g$ to the hyperplane 
$\{(x_1,x_2,x_3): x_3 = y_3\}$.  Clearly
$$ 
\sum_{y_3 \in F} \| g_{y_3} \|_{L^{8/5}(F^3, dx)}
\leq |F|^{3/8}\ (\sum_{y_n \in F} \| g_{y_3} \|_{L^{8/5}(F^3,dx)}^{8/5}\ )^{5/8}
= |F|^{3/8}\ \| g \|_{L^{8/5}(F^3, dx)}
$$
so by the triangle inequality it suffices to show that
\be{g-targ}
\| g_{y_3} * K \|_{L^4(F^3, dx)} \lessapprox |F|^{-1/4}\ \| g_{y_3} \|_{L^{8/5}(F^3, dx)}
\end{equation}
for each $y_3$.  By translation invariance we may take $y_3 = 0$.

From \eqref{gauss} we have
$$ 
K(\underline{x},x_3) = |F|^{-2}\ S(x_3)^{2}\
e(\underline{x} \cdot \underline{x} / 4 x_3)
$$ 
so we can write the left-hand side of \eqref{g-targ} as
$$ 
\Big(\sum_{x_3 \in F: x_3 \neq 0} \sum_{\underline{x} \in F^2}
|\, |F|^{-2}\ |S(x_3)|^{2}
\sum_{\underline{y} \in F^2} g(\underline{y},0)\
e((\underline{x} - \underline{y})\cdot (\underline{x} - \underline{y}) / 4 x_3)
\ |^4\,\Big)^{1/4}.$$
By \eqref{gauss-size} we have $|S(x_3)|^2 = |F|$ for $x_3\neq 0$.  We now make the 
``pseudo-conformal'' substitution $t := 1/4x_3$ and 
$z := -\underline{x}/2x_3$, so that
$$ (\underline{x} - \underline{y}) \cdot (\underline{x} - \underline{y}) / 4 x_3
= z^2 x_3 + z \cdot \underline{y} + t \underline{y}^2$$
and the previous expression becomes
$$ 
|F|^{-1}\,  
\Big(\sum_{t \in F: t \neq 0} \sum_{z \in F^2}
|\ e(z^2 x_3)
\sum_{\underline{y} \in F^2} g(\underline{y},0)\ \,
e((z,t) \cdot (\underline{y}, \underline{y}^2))\
|^4\,\Big)^{1/4}.
$$
The phase $e(z^2 x_3)$ can be discarded, as can the restriction 
$t \neq 0$.  We can thus bound the previous by
$$
|F|^{-1} \big(\sum_{(z,t) \in F^3} | (G d\sigma)\spcheck(z,t)|^4\, \big)^{1/4}
$$
where
$$ 
G(\underline{y}, \underline{y}^2) := |F|^2\   g(\underline{y},0).
$$
By Proposition \ref{12-7} we may bound this by
$$ 
\lessapprox |F|^{-1}\ \| G \|_{L^{8/5}(S, d\sigma)}
= |F|^{-1/4}\ \| g_0 \|_{L^{8/5}(F^3, dx)},
$$
and the claim follows.
\end{proof}
The logarithmic loss can probably be removed.
Since the above restriction theorem uses one restriction estimate to 
prove another, it may be possible to iterate it (as in 
\cite{tvv:bilinear}, \cite{tv:cone1}) to obtain some improvement, but 
this is unlikely to reach the conjectured best possible bound of 
$R^*(2 \to 3)$.  This estimate also compares favorably with the best 
bound for the paraboloid in $\R^3$, which is $R^*(p \to 26/7+)$ for 
certain $2 < p \leq \infty$, and it may be that the argument above 
also has some application to the Euclidean problem.

The question of what $L^p \to L^q$ estimates the convolution operator 
$f \mapsto f * K$ is itself quite interesting, being the finite field 
version of the Bochner-Riesz problem (see e.g. \cite{borg:kakeya}, 
\cite{stein:large} for a discussion).  From the above argument we see 
that this problem is closely related to restriction, especially for 
paraboloids.  (In the Euclidean case there are further connections, 
see \cite{tao:boch-rest}, although the arguments there rely on scaling 
and so do not extend to the finite field case).  We will not discuss 
this topic further here though.

\section{The cone}\label{cone-sec}

We now consider the restriction problem for the cone
\be{cone}
S := \{ (\xi, u, v): \xi, u, v \in F_*; uv = \xi^2 \} \setminus \{ (0,0,0) \}
\end{equation}
in $F^3$; the higher-dimensional cones are also of interest but will not be 
discussed here.  We have removed the origin $(0,0,0)$ for technical 
convenience, but it can be restored with no significant change to the 
results.

It is relatively easy to obtain the boundedness of $R^*(2 \to 4)$ (the 
Euclidean counterpart of this is in \cite{barcelo}, although the exponents 
are slightly different).  One expects\footnote{The decay here is 
$\tilde d=1$ and so a Tomas-Stein style argument would only yield the 
boundedness of $R^*(2 \to 6)$.} to apply Lemma \ref{k-planch}, however 
there is a slight difficulty because the cone contains lines through 
the origin, which is bad for \eqref{sum} when $\eta = 0$.  However, 
when $\eta \neq 0$ a routine algebraic computation shows that the number 
of solutions to \eqref{sum} is $O( |F| )$ (in other words, the 
intersection of a cone with a non-trivial translate of itself is at 
most one-dimensional).  By the argument in the proof of 
Lemma \ref{k-planch} we thus have
$$ 
\| f d\sigma * f d\sigma \|_{L^{2}(F_*^3 \backslash \{0\}, d\xi)} \lesssim  
\| f \|_{L^2(S, d\sigma)}^2.
$$
On the other hand, a direct computation shows
$$ 
|fd\sigma * fd\sigma(0)| \lesssim |F| \| f \|_{L^2(S, d\sigma)}^2,
$$
and the claim follows by summing.

Because the cone contains one-dimensional spaces (if we re-insert 
the origin) we obtain from \eqref{q-improv} the necessary 
condition $q \geq 2p'$ for the boundedness of $R^*(p \to q)$.  
Thus the bound on $R^*(2 \to 4)$ is sharp in one sense.  In fact 
it is the best restriction bound possible:

\begin{proposition}\label{counter}  If $q < 4$ and $1 \leq p \leq \infty$, 
then $R^*(p \to q)$ is unbounded.
\end{proposition}

\begin{proof}
The idea is to complete the square in the phase of the cone restriction 
operator and then use the standard estimate for Gauss sums.

We apply \eqref{rest-dual} with $f := \chi_X$, where $X$ is the set
$$ 
X := \{ (x,y,z) \in F^3: z \hbox{ is a non-zero square, and } y = x^2 / 4z \}.
$$
Clearly $|X| \sim |F|^2$.  We will show that
\be{hat}
\| \hat \chi_X \|_{L^{p'}(S, d\sigma)} \gtrsim |F|^{3/2},
\end{equation}
which implies the claim since $q > 4$.  In fact we will show that
$$ |\hat \chi_X(\xi, u, v)| \sim |F|^{3/2}$$
for all $(\xi, u, v) \in S$ with $u \neq 0$.

To see this, we write $v = \xi^2/u$ and $\xi=tu$, 
let $Q$ be the set of all non-zero 
squares in $F$, and compute
\bas
\hat \chi_X(tu, u, t^2u) &=
\sum_{z \in Q} \sum_{x \in F} e(x tu + x^2 u/4z + t^2 u z)\\
&=
\sum_{z \in Q} \sum_{x \in F} e(u(x + 2 t z)^2/4z)\\
&=
\sum_{z \in Q} \sum_{y \in F} e(u y^2) \quad \hbox{(substituting } y := (x + 2tz)/\sqrt{2z})\\
&= |Q|\ S(u)
\end{align*}
and the claim follows from \eqref{gauss-size} and $2|Q|=|F|-1$.
\end{proof}

Unlike the situation in Euclidean space $\R^3$, this counter-example does 
not seem to be removable simply by passing to a bilinear formulation.  
However, it may well be that the large body of work on the Euclidean cone 
restriction problem (\cite{borg:cone}, \cite{mock:cone}, 
\cite{wolff:cone}, \cite{tv:cone1}, etc.) still has some application to 
the finite field case.

We summarize our results on the restriction problem for finite fields in 
Figure \ref{fig1}.

\begin{figure}\label{fig1}
\begin{tabular}{|l|l|l|} \hline
Surface & Best known  			& Best known \\
{}	& restriction theorem            & counterexample
\\\hline

$(t, \ldots, t^n)$, $\chr (F) > n$ & $R^*(2 \to 2n)$ & $R^*(2 \to 2n)$ \\
$n=2,$ parabola & $R^*(2 \to 4)$ & $R^*(2 \to 4)$ \\
$n=3,$ paraboloid, $-1$ non-square & $R^*(8/5 \to 4)$, $R^*(2 \to 18/5+)$ & $R^*(2 \to 3)$ \\
$n=3,$ paraboloid, $-1$ square & $R^*(2 \to 4)$ & $R^*(3 \to 3)$ \\
$n=3,$ cone & $R^*(2 \to 4)$ & $R^*(2 \to 4)$ \\
\hline
\end{tabular}
\caption{Some of the surfaces discussed in this paper, the best restriction 
theorem we could prove (up to logarithms), and the conjectured best 
restriction theorem suggested by the counterexamples. }
\end{figure}

\section{The Kakeya problem}\label{kakeya-sec}

Having completed our discussion of restriction problems for the moment, 
we now set up the notation for the Kakeya problem.

We parameterize the vector space $F^n$ by $x = (\underline{x},x_n)$, where 
$\underline{x} \in F^{n-1}$, $x_n \in F$.  For any $x_0, v \in F^{n-1}$, 
define the line $l(x_0,v)$ by 
$$ l(x_0,v) := \{ (x_0 + vt,t): t \in F \}.$$
We refer to $v$ as \emph{directions}, and endow the space $F^{n-1}$ of 
directions with normalized counting measure
$$ 
\int_{F^{n-1}} f(v)\ dv := \frac 1{|F|^{n-1}} \sum_{v \in F^{n-1}} f(v).
$$
In this definition we have excluded the horizontal lines, but this will 
make no essential difference to our results.

Define a \emph{Besicovitch set} to be any subset $E$ of $F^n$ which 
contains a line in every direction, i.e. for every $v \in F^{n-1}$ there 
exists an $x_0$ such that $l(x_0,v) \in E$.

As an example of a Besicovitch set, and consider 
\be{2-d}
E := \{ (x,t) \in F^2: x + t^2 \hbox{ is a square } \}.
\end{equation}
This set has cardinality $\frac{1}{2} |F|^2 + \frac{1}{2} |F|$, and for 
every $v \in F$ the line $l(v^2/4, v)$ is contained in $E$.  (cf. the 
``completing the square'' trick in Proposition \ref{counter}).

The \emph{Kakeya set conjecture} for finite fields asserts that every 
Besicovitch set has cardinality $|E| \approx |F|^n$.  Informally, this 
means that it is impossible to compress lines in distinct directions 
into a small set.

This conjecture is proven in two dimensions but is open in higher 
dimensions.  In three dimensions the best bound is $|E| \gtrsim |F|^{5/2}$ 
(see \cite{wolff:survey} or the arguments below).  For higher dimensions, 
see below.

There are many important variations of the Kakeya problem, in which lines 
are replaced by circles, planes, light rays, spheres, or other geometric 
objects, or if the requirement of distinct directions is replaced by some 
other condition.  We do not attempt a survey of all the possibilities here, 
but refer the reader to \cite{wolff:survey}.

One obvious attempt to construct a counterexample to the Kakeya set 
conjecture would be to set $E$ equal to some algebraic variety such 
as $\{ x \in F^n: P(x) = 0 \}$ where $P: F^n \to F$ is some polynomial.  
However, this cannot work:

\begin{proposition}\label{algebraic}  Let $K > 0$, and let $E$ be 
contained in the set
$$ \bigcup_{i=1}^K \{ x \in F^n: P_i(x) = 0\},$$
where for each $1 \leq i \leq K$, $P_i: F^n \to F$ is a non-zero polynomial 
of degree at most $K$.  Then, if $\chr(F)$ is sufficiently large depending 
on $K$, the set $E$ cannot be a Besicovitch set.
\end{proposition}

\begin{proof}
Suppose for contradiction that $E$ is a Besicovitch set.  Then by the 
pigeonhole principle for every $v \in F^{n-1}$, there exists $1 \leq i 
\leq K$ such that $l(x_0(v),v)$ intersects $\{ x \in F^n: P_i(x) = 0 \}$ 
in at least $|F|/K$ points.  By another pigeonholing, we may therefore 
find $1 \leq i \leq K$ and a set $V \subset F^{n-1}$ of cardinality 
$|V| \geq |F|^{n-1}/K$ such that
$$ |l(x_0(v),v) \cap \{ x \in F^n: P_i(x) = 0 \}| \geq |F|/K$$
for all $v \in V$.  In other words,
$$ |\{ t \in F: P_i(x_0(v) + tv, t) = 0 \}| \geq |F|/K \hbox{ for all } v \in V.$$
If $\chr(F)$ is sufficiently large, this implies (since $P_i$ has bounded 
degree) that
$$ P_i(x_0(v) + tv, t) = 0 \hbox{ for all } v \in V \hbox{ and } t \in F.$$
Suppose that $P_i$ has degree $d$ for some $0 < d \leq M$, and let $P_i^*$ 
be the principal part of $P_i$ (i.e. the terms which have degree exactly 
equal to $d$).  Then $P_i(x_0(v) + tv, t)$ is a polynomial in $t$ of 
degree $d$.  Extracting the $t^d$ component we see that
$$ P_i^{*}(tv, t) = 0 \hbox{ for all } v \in V \hbox{ and } t \in F.$$
Fix $t$, and think of $P_i^{*}(tv, t)$ as a polynomial in $v$ of degree 
at most $d$.  This is zero on at least $1/K$ of the values of $F^n$, which 
implies (if $\chr(F)$ is sufficiently large) that $P_i^*(tv,t)$ is 
identically zero.  But this implies that $P_i^*(x,t) = 0$ for all 
$(x,t) \in F^n$, which is absurd since $P_i$ was supposed to have degree 
$d$.  This is the desired contradiction.
\end{proof}

Other attempts to produce counterexamples, e.g. by using the squares as 
in \eqref{2-d}, or by random constructions, do not appear to give 
Besicovitch sets which are much smaller than $F^n$.  However this is far 
from a proof that the Kakeya set conjecture is true.

The Kakeya set conjecture can be attacked directly, however we shall 
study it in the context of a more general conjecture concerning maximal 
functions.

For any function $f$ on $F^n$, define the \emph{Kakeya maximal function} 
$f^*$ on $F^{n-1}$ by
\be{max-def}
f^*(v) = \sup_{x_0 \in F^{n-1}} \sum_{x \in l(x_0,v)} |f(x)|.
\end{equation}

Let $1 \leq p,q \leq \infty$ be exponents.  We define $K(p \to q)$ to be 
the best constant such that 
\be{kak-est}
\| f^* \|_{L^q(F^{n-1}, dv)} \leq K(p \to q)\ \| f\|_{L^p(F^n, dx)}
\end{equation}
for all $f$.  By linearization and duality, this estimate is equivalent to 
the statement that
\be{kakeya-alt}
\| \int_{v \in F^{n-1}} g(v)\  \chi_{l(x_0(v),v)}\ dv \|_{L^{p'}(F^n, dx)}
\leq K(p \to q)\ \| g \|_{L^{q'}(F^{n-1}, dv)}
\end{equation}
for all $g$ on $F^{n-1}$, and all functions $x_0: F^{n-1} \to F^{n-1}$.  
Clearly we may take $f$ and $g$ non-negative in the above.

For instance, it is easy to verify that
$$ 
K(1 \to q) = 1 \hbox{ and } K(\infty \to q) = |F|
$$
for all $1 \leq q \leq \infty$, while from the H\"older and Young 
inequalities we see that $K(p \to q)$ is non-decreasing in both $p$ 
and $q$.  By testing \eqref{kak-est} with $f$ equal to the characteristic 
function of a point, line, or all of $F^n$, we obtain the bounds
$$ 
K(p \to q) \geq 1, |F|^{-(n-1)/q + 1/p'}, |F|^{1-n/p}.
$$
Thus in order for $K(p \to q)$ to be bounded one must have
$$ 
p \leq n \hbox{ and } q \geq (n-1)p'.
$$
The \emph{Kakeya maximal conjecture} for finite fields asserts that these 
necessary conditions are also sufficient.  In particular, $K(n \to n)$ 
should be bounded.  The continuous version of this conjecture has been 
intensively studied; see \cite{wolff:survey} for a survey.

This conjecture is related to the Kakeya set conjecture in the following 
sense: if $K(p \to q)$ is bounded, then Besicovitch sets have cardinality 
at least $\gtrsim |F|^p$.  To see this, observe from construction that if 
$E$ is a Besicovitch set, then
$$ (\chi_E)^*(v) = |F|$$
for all directions $v$.  Inserting this into \eqref{kak-est} the claim 
follows.

Thus it is of interest to bound $K(p \to q)$ for as large a value of $p$ 
as possible.  When $p=2$ this is quite easy:

\begin{proposition}[C\'ordoba's argument]\label{cordoba} $K(2 \to 2n-2) 
\leq \sqrt{2}$.
\end{proposition}

\begin{proof}
We use \eqref{kakeya-alt}, and compute
$$
\| \int_{v \in F^{n-1}} g(v)\ \chi_{l(x_0(v),v)}\ dv \|_{L^2(F^n, dx)}^2
$$
as
$$
\int_{v,v' \in F^{n-1}} g(v)\ g(v')\ |l(x_0(v),v) \cap l(x_0(v'),v)|\ dv dv'.$$
First consider the diagonal contribution $v=v'$.  This is
$$
\frac 1{|F|^{n-1}} \int_{F^{n-1}} g(v)^2\ |F|\ dv
= |F|^{2-n}\ \| g \|_{L^2(F^{n-1}, dv)}^2
\leq \| g \|_{L^{(2n-2)'}(F^{n-1}, dv)}^2.
$$
Now consider the off-diagonal contribution $v \neq v'$.  Since two 
non-parallel lines can only intersect in at most one point, this 
contribution is at most
$$
\int_{v,v' \in F^{n-1}} g(v)\ g(v') \ dv dv' = \| g \|_{L^1(F^{n-1}, dv)}^2
\leq \| g \|_{L^{(2n-2)'}(F^{n-1}, dv)}^2.$$
Adding the two terms, the claim follows.
\end{proof}

Observe that this estimate is on the line $q = (n-1)p'$, and implies that 
Besicovitch sets have cardinality at least $\frac{1}{2} |F|^2$ (compare 
with the example \eqref{2-d}).  This settles the Kakeya conjectures in 
two dimensions.

In the remainder of our discussion we shall be prepared to lose logarithmic 
factors.  In this case the Kakeya problem is equivalent to the classic 
problem of counting incidences between points and lines, but with the 
condition that the lines all point in different directions.

\begin{proposition}\label{equiv}  Let $1 \leq p,q \leq \infty$.  The 
statement $K(p \to q) \lessapprox 1$ holds if and only if one has the 
incidence bound
\be{incidence}
|\{ (p,l) \in P \times L: p \in l \}| \lessapprox |P|^{1/p}\ |L|^{1/q'}\ 
|F|^{(n-1)/q}
\end{equation}
for all collections $P$ of points in $F^n$ and all collections $L$ of 
lines in $F^n$, each of which points in a different direction.
\end{proposition}

\begin{proof}
Suppose first that $K(p \to q) \lessapprox 1$, and let $P$ and $L$ be as 
above.  Let $V \subset F^{n-1}$ denote the directions of the lines in $L$, 
then we see that
$$ |\{ (p,l) \in P \times L: p \in l \}| \leq \sum_{v \in V} 
(\chi_P)^*(v) = |F|^{n-1} \int_V (\chi_P)^*.$$
Applying H\"older and \eqref{kak-est} and noting that $|V| = |L|$, the 
claim \eqref{incidence} follows.

Now suppose that \eqref{incidence} holds.  By \eqref{kakeya-alt} and 
duality it suffices to show that
$$
\sum_{x \in F^n} \int_{v \in F^{n-1}} g(v)\ \chi_{l(x_0(v),v)}(x)\ f(x)\ dv 
\lessapprox \| g \|_{L^{q'}(F^{n-1}, dv)}\ \| f \|_{L^p(F^n, dx)}.
$$
By the standard dyadic pigeonholing argument we may assume that $f = \chi_P$ 
and $g = \chi_V$ are characteristic functions.  (Note that since there are 
at most $|F|^n$ points and directions, we only need to pigeonhole into 
$O(\log |F|)$ categories).  The above estimate then becomes
$$ 
|F|^{1-n} \sum_{x \in P} \sum_{v \in V}\ \chi_{l(x_0(v),v)}(x)
\lessapprox |F|^{(1-n)/q'}\ |V|^{1/q'}\ |P|^{1/p}.
$$
But this follows from \eqref{incidence} by setting 
$L := \{ l(x_0(v),v): v \in V \}$.
\end{proof}

As discussed before we have the constraint $q \geq (n-1)p'$ as a necessary 
condition for \eqref{incidence} to hold.  However if one is prepared to lose 
the trivial error term of $|P| + |L|$ then one can do much better; see 
Proposition \ref{lines}, or Proposition \ref{wolff} below.  To prepare 
for this we give

\begin{corollary}\label{implic}  Suppose that we can prove a bound of the form
\be{incid-2}
|\{ (p,l) \in P \times L: p \in l \}| \lessapprox |P|^{a}\ |L|^{1-b}\ |F|^{1-c} + |P| + |L|
\end{equation}
for some $0 \leq a,b,c \leq 1$ with $(n-1)b+c \geq 1$, whenever $L$ has 
distinct directions.  Then we have
$$ K(p \to q) \lessapprox 1$$
where $p := \frac{1}{a}((n-1)b + c)$ and $q := \min((n-1)p', \frac{1}{b}(n-1)b+c))$.
\end{corollary}

\begin{proof}  Denote the left-hand side of \eqref{incid-2} by $|I|$.
By \eqref{incidence} it suffices to prove
\be{i-targ}
\frac{|I|}{|L|\, |F|} \lessapprox |P|^{1/p}\ \, 
\Big(\frac{|L|}{|F|^{n-1}}\Big)^{1/q}\ \,|F|^{-1}.
\end{equation}
If $|I| \lesssim |P|$ then \eqref{i-targ} follows from the trivial bounds 
$|P| \leq |F|^n$ and $|L| \geq 1$ since $q \geq (n-1)p'$.  If 
$|I| \lesssim |L|$ then \eqref{i-targ} follows from the trivial bounds 
$|L|/|F|^{n-1} \leq 1$ and $|P| \geq 1$.  By hypothesis we may therefore 
assume
$$ 
\frac{|I|}{|L|\, |F|} \lessapprox 
|P|^a\ \,  \Big(\frac{|L|}{|F|^{n-1}}\Big)^{-b}\ \, |F|^{-(n-1)b-c}.
$$
The claim then follows by interpolating this with the trivial bound 
$|I| \leq |L| |F|$ (since every line contains at most $|F|$ points).
\end{proof}

Thus, for instance, Proposition \ref{lines} can be used to imply 
Proposition \ref{cordoba} in the two dimensional case (at least if 
one is willing to lose logarithms).  It is an interesting question 
to ask whether Proposition \ref{lines} can be improved.  Consider 
for instance the case $|L| = |P|$, in which case the Proposition 
gives an incidence bound of $|P|^{3/2}$.  This is sharp when $|P| = |F|^2$ 
(just set $P$ and $L$ equal to the space of all points and all lines), and 
more generally for any subfield $G$ of $F$ the estimate is sharp when 
$|P| = |G|^2$ (set $P$ equal to $G^2$, and $L$ equal to lines with 
slope and intercept in $G$).  However, one expects to do better when 
$|F|$ does not contain any sub-fields.  For instance, if $F = \Z/p\Z$ 
for some prime $p$ and $|P| = |L| = p$, is it possible to substantially 
improve the incidence bound of $p^{3/2}$?  In the Euclidean case one 
would obtain a bound of $O(p^{4/3})$ from the famous Szemer\'edi-Trotter 
theorem \cite{szemeredi:trotter}, but this argument uses crucially the 
fact that $\R$ is an ordered field and so does not apply to $\Z/p\Z$.  

It turns out that the question of improving the bound of $p^{3/2}$ is 
equivalent to disproving the existence of a set $E \subset \Z/p\Z$ for 
which $|E| \approx |E+E| \approx |E \cdot E| \approx \sqrt{p}$.  
See \cite{katz:falc}, \cite{katz:falc-expos}, \cite{elekes}.  In other words, the question is 
equivalent to whether $\Z/p\Z$ contains an ``approximate sub-ring'' of 
``dimension $1/2$''.

There are similar connections between the continuous versions of these 
problems, specifically the Falconer distance problem, the Furstenburg 
set problem, and the Erd\"os ring problem.  See \cite{katz:falc}, \cite{katz:falc-expos}.

The argument in Proposition \ref{lines} extends without difficulty to 
higher dimensions\footnote{Of course, one loses duality of points and 
lines in higher dimensions, however the only geometric fact we needed 
was that every two distinct lines are incident to at most one point, 
and this fact is preserved on the dual side (i.e. two distinct points 
determine at most one line).}, and one can thus conclude using 
Proposition \ref{equiv} that $K((n+1)/2 \to (n+1)) \lessapprox 1$.  
This is the analogue of the results in \cite{cdrdef} (see also 
\cite{drury:xray}), and shows that Besicovitch sets have cardinality 
$\gtrapprox |F|^{(n+1)/2}$.

We now improve this to $(n+2)/2$ by the following argument of Wolff 
(\cite{wolff:kakeya}, \cite{wolff:survey}; we shall use an argument 
due to Nets Katz).

\begin{definition}\label{lax-def} A collection $L$ of lines in $F^n$ 
obeys the \emph{Wolff axiom} if every $2$-plane contains at most $O(|F|)$ 
lines from $L$.
\end{definition}

Clearly any collection $L$ of lines with distinct directions will obey the 
Wolff axiom.

\begin{proposition}\label{wolff}  If $L$ obeys the Wolff axiom, then
\be{Wolff-est}
|\{ (p,l) \in P \times L: p \in l \}| \lesssim |P|^{1/2}\ |L|^{3/4} \
|F|^{1/4} + |P| + |L|.
\end{equation}
\end{proposition}

\begin{proof}
Denote the set on the left-hand side of \eqref{Wolff-est} by $I$.  We may 
assume without loss of generality that we have the ``two-ends condition''
\be{large}
|I| \gg |P| 
\end{equation}
(i.e. most lines contain at least two points) and the ``bilinear condition''
\be{bilinear}
|I| \gg |L| 
\end{equation}
(i.e. most points lie in at least two lines) since \eqref{Wolff-est} is 
trivial otherwise.

The idea will be to obtain a lower and upper bound for the number of 
triangles in $F^n$ formed by $P$ and $L$.

Let $L'$ denote the set of lines with the average or above-average number of 
points in $P$:
$$ 
L' := \{ l \in L: |l \cap P| \geq |I|/2\ |L| \}.
$$
Let
$I' := \{ (p,l) \in I: l \in L' \}$.
Clearly we have
$| I \backslash I'| < |I|/2$,
so that
\be{isimi}
|I'| \sim |I|.
\end{equation}

Let $V$ denote the set of (possibly degenerate) ``angles''
$$ V := \{ (p,l,l') \in P \times L' \times L': (p,l),(p,l') \in I' \}.$$
From \eqref{cauchy} we have 
$$ |V| \geq \frac{{|I'|}^2}{|P|} \sim \frac{|I|^2}{|P|}.$$
Let $V'$ denote the non-degenerate angles
$V' := \{ (p,l,l') \in V: l \neq l' \}$.
Clearly
$|V \backslash V'| = |I|$.
By \eqref{large} and the preceding we thus have
\be{v-est}
|V'| \gtrsim \frac{|I|^2}{|P|}.
\end{equation}
Let $W$ denote the set of non-degenerate pointed angles
$$ 
W = \{ (p,l,l',p') \in V' \times P: p' \in l', p \neq p' \}.
$$
From the definition of $L'$, each element of $V'$ contributes at least 
$\frac{|I|}{2|L|} - 1$ elements to $W$.  From \eqref{bilinear} we thus have
\be{w-est}
|W|\gtrsim \frac{|I|\, |V'|}{|L|}.
\end{equation}
Let $T$ denote the space of pairs of linked non-degenerate pointed angles:
$$ 
T = \{ (p_1,l_1,l'_1,p'_1, p_2, l_2, l'_2,p'_2) \in W \times W: 
l_1 = l_2, p'_1 = p'_2 \}.
$$
The set $T$ can be identified with the space of triangles with points and 
lines in $P$ and $L$, although there is a possible degeneracy in that $p_1$ 
may equal $p_2$.  From \eqref{w-est} and \eqref{cauchy} we thus have
$$
|T| \geq \frac{|W|^2}{|P||L|}.
$$
Now let $T'$ be the space of non-degenerate triangles:
$$ 
T' = \{ (p_1,l_1,l'_1,p'_1, p_2, l_2, l'_2,p'_2) \in T: p_1 \neq p_2 \}.
$$
The set $T \backslash T'$ is bijective with $W$.  From the previous 
estimates and \eqref{large} we thus have 
$|T \backslash T'| = |W| \ll |T|$, thus
$$ 
|T'| \gtrsim \frac{|W|^2}{|P|\,|L|}.
$$
This is our lower bound for the number of triangles.  To obtain an upper 
bound,
we observe from the Wolff axiom that $p_1, l_1, l'_1$ contributes at most 
$|F|$ elements to $T'$.  This is because the line $l'_2$ lies in the plane 
generated by $l_1$, $l'_1$.  The remaining variables are determined by
 $\{ p'_1\} = \{ p'_2 \} = l'_1 \cap l'_2$, $\{ p_2\} = l_2 \cap l'_2$, 
and $l_2 = l_1$.  Since $(p_1,l_1,l'_1) \in V'$, we thus have
$$ 
| T'| \leq |V'| |F|.
$$
Combining the previous two estimates with \eqref{w-est} we obtain
$$ 
|V'| |F| \gtrsim \frac{|I|^2 |V'|^2}{|P| |L|^3}.
$$
Dividing by $|V'|$ and then applying \eqref{v-est} we obtain
$$  
|F| \gtrsim \frac{|I|^4}{|P|^2 |L|^3}
$$
and \eqref{Wolff-est} follows.
\end{proof}

From the above and Corollary \ref{implic} we thus have\footnote{It is clear 
from the proof of Corollary \ref{implic} that \eqref{Wolff-est} in fact is 
substantially stronger than this Kakeya estimate, especially with the 
two-ends and bilinear assumptions.  This is consistent with experience 
in the Euclidean case, see e.g. \cite{wolff:kakeya} for a discussion of 
the two-ends property, and \cite{tvv:bilinear} for the bilinear property.} 
that 
\be{wk}
K\big(\frac{n+2}{2} \to \frac{(n-1)(n+2)}{n}\big) \lessapprox 1.
\end{equation}
This is the analogue of the main result of \cite{wolff:kakeya}, and also 
lies on the line $q = (n-1)p'$.  In particular, Besicovitch sets have 
cardinality $\gtrapprox |F|^{(n+2)/2}$.  In fact one can remove the 
logarithm by appealing directly to Proposition \ref{wolff}.

In the important three-dimensional case we do not know if this bound of 
$(n+2)/2 = 5/2$ can be improved.  However, we cannot improve this bound 
solely by using the Wolff axiom and the cardinality bound on $L$. To see 
this, suppose that $F$ has a subfield $G$ of index 2, and let $P$ be the 
\emph{Heisenberg group}
$$ P := \{ (z_1, z_2, z_3) \in F: \Im(z_1 \overline{z_2}) = \Im(z_3) \}$$
where $z \mapsto \overline{z}$ is the non-identity involution which 
preserves $G$, and $\Im(z) := (z - \overline{z})/2$.  The set $P$ has 
cardinality $\sim |F|^{5/2}$, and contains the family of lines
$$ 
L := \{ l((x_1,x_2),(v_1,v_2)): \Im(v_1 \overline{v_2}) = \Im(x_1 
\overline{x_2}) = 0; v_1 \overline{x_2} + v_2 \overline{x_1} = 1 \}.
$$
This family of lines can be seen to have cardinality $\sim |F|^2$ and 
obeys the Wolff axiom, but the lines do not all point in distinct 
directions.  Thus to improve the $(n+2)/2$ bound in three dimensions 
one must somehow use the distinctness of the directions more 
intimately\footnote{Alternatively, if one is working in a field such 
as $F = \Z/p\Z$, one could try to preclude approximate half-dimensional 
rings of the type mentioned earlier, in order to eliminate Heisenberg-type 
counterexamples.}. On the other hand, it may be possible to improve 
$(n+2)/2$ in higher dimensions just by using the Wolff axiom (together 
with the obvious generalizations to subspaces of dimension greater than 2).

One possible approach to the Kakeya problem which does indeed exploit 
distinctness of directions is via arithmetic combinatorics (see
 \cite{borg:high-dim}, \cite{katz:high-dim}, \cite{kt:advanced}).  We 
borrow the following notation from \cite{kt:advanced}.

Define a \emph{slope} $r$ to be any element of 
$$ 
\{ \frac{a}{b}: a,b \in \Z: 0 < b < \chr(F) \} \cup \{\infty\}.
$$
We define the projections $\pi_r: F^{n-1} \times F^{n-1} \to F^{n-1}$ by
$\pi_r(a,b) := a + rb$ for $r \neq \infty$ and $\pi_\infty(a,b) := b$.
We call a slope \emph{proper} if $r \neq -1$.

\begin{definition}\label{sd-def}  Let $R$ be a finite collection of proper 
slopes, and let $\alpha \in \R$.  We say that the statement $SD(R, \alpha)$ 
holds if one has the bound
$|G| \lesssim \sup_{r \in R} |\pi_r(G)|^\alpha$
whenever $G \subseteq F^{n-1} \times F^{n-1}$ is a finite set obeying 
\be{one-to-one}
\pi_{-1} \hbox{ is one-to-one on } G.
\end{equation}
\end{definition}

One can think of the statement $SD(R,\alpha)$ as a quantitative way to 
control the cardinality of the projection $\{ a-b: (a,b) \in G\}$ in 
terms of other projections $\{ a+rb: (a,b) \in G\}$.

One trivially has $SD(R, 2)$ as soon as $R$ contains at least two elements.  
In \cite{borg:high-dim} Bourgain improved this to $SD(\{0, 1, \infty\}, 
2 - \frac{1}{13})$ based on some arguments of Gowers \cite{gowers}; this 
was later improved to $SD(\{0, 1, \infty\}, 2 - \frac{1}{6})$ in 
\cite{katz:high-dim}, or $SD(\{0, 1, 2, \infty\}, 2 - \frac{1}{4})$ 
if $\chr(F)$ is sufficiently large.  The current record is that for 
any $F$ there exists a finite set $R$ and an $\eps_F > 0$ such that 
$\eps_F \to 0$ as $\chr(F) \to \infty$ and $SD(R, \alpha + \eps_F)$ 
holds, where
$\alpha = 1.675\ldots$ is the largest root of $\alpha^3 - 4\alpha + 2=0$ 
(see \cite{kt:advanced}; the statement there is for vector spaces but 
the extension to finite fields of sufficiently large characteristic is 
routine).

These problems are related to the Balog-Szemeredi theorem \cite{balog} 
(see also \cite{gowers}).  Their connection to Kakeya lies through 
relationships such as the following.

\begin{lemma}\label{slices}  If $SD(R,\alpha)$ holds for some $R$, 
then Besicovitch sets have cardinality $\gtrsim |F|^{(n-1)/\alpha + 1}$ 
(the implicit constant depends on $R$).
\end{lemma}

In particular, if one could prove sums-differences lemmas with $\alpha$ 
arbitrarily close to one, one would be able to settle the Kakeya set 
conjecture.  To improve upon the three-dimensional results one would 
need $\alpha$ to be better than $4/3$.

\begin{proof}  We follow the arguments of Bourgain \cite{borg:high-dim}.  
Let $E$ be a Besicovitch set, and let $C_0$ be a large number (depending 
on $R$) to be chosen later.  Call a height $t \in F$ \emph{exceptional} if
$$ 
| \{ x \in F^{n-1}: (x,t) \in E \}| \geq C_0 |E|/|F|.
$$
Clearly at most $|F|/C_0$ heights are exceptional.  By a simple 
probabilistic argument we thus see that if $C_0$ is sufficiently large, 
then we may find distinct heights $t_0, t_\infty \in F$ such that the heights
$$ 
t_r := \frac{1}{r+1} t_0 + \frac{r}{r+1} t_\infty
$$
are non-exceptional for all $r \in R$.

Fix $t_0, t_\infty$, and define
$$ 
G := \{ (x_0(v) + t_0 v, x_0(v) + t_\infty v): v \in F^{n-1} \}.
$$
Since $t_0 \neq t_\infty$, we see that $\pi_{-1}$ is one-to-one on $G$ and
$$ 
|G| = |F|^{n-1}.
$$
On the other hand, since $(x_0(v) + t_r v, t_r) \in E$ for all 
$v \in F^{n-1}$ and $r \in R$, we have
$$ 
|\pi_r(G)| \leq|\{ x \in F^{n-1}: (x,t_r) \in E \}| \leq C_0 |E|/|F|.
$$
Applying $SD(R,\alpha)$ the claim follows.
\end{proof}

This Lemma immediately gives some bounds on the size of Besicovitch sets, 
for instance the result $SD(\{0, 1, \infty\}, 2-\frac{1}{4})$ gives 
$|E| \gtrsim |F|^{(4n+3)/7}$ for sufficiently large values of $\chr(F)$.  
We give a direct proof of this fact below:

\begin{proposition}\label{47} \cite{katz:high-dim} If $\chr(F) > 3$, then 
we have $|E| \gtrsim |F|^{(4n+3)/7}$.
\end{proposition}

\begin{proof}
We set $P$ equal to $E$, and $L$ equal to the lines 
$L = \{ l(x_0(v),v): v \in F^{n-1} \}$.  Assume for contradiction 
that $|P| \ll |F|^{(4n+3)/7}$.

Let $I = \{ (p,l) \in P \times L: p \in L \}$ be the set of incidences, 
then $|I| = |F|^n$.  Let $V'$ denote the non-degenerate angles
$$ 
V' := \{ (p,l_1,l_2) \in P \times L \times L: p \in l_1; p \in l_2; 
l_1 \neq l_2 \}.
$$
Arguing as in Proposition \ref{wolff} we have
$$
|V'| \gtrsim \frac{|I|^2}{|P|} = |F|^{2n}\ |P|^{-1}.
$$
Thus if we let $A$ denote the set
$$ 
A := \{ (p,p_1,p_2,l_1,l_2) \in P^3 \times L^2: p, p_1 \in l_1; p, p_2 \in l_2; l_1 \neq l_2; p \neq p_1; p \neq p_2 \}
$$
we have
$$ | A| \gtrsim |F|^2\ |V'| \gtrsim |F|^{2n+2}\ |P|^{-1}.$$
The pair $(p_1,p_2)$ lies in a set of cardinality $|P|^2$, so by 
\eqref{cauchy} we have
$$ 
| Q |\gtrsim \frac{|A|^2}{|P|^2} \gtrsim |F|^{4n+4} |P|^{-4}
$$
where $Q$ denotes the set of quadrilaterals
\bas
 Q := \{ &(p,p',p_1,p_2, l_1,l_2,l'_1,l'_2) \in P^4 \times L^4:\\
&p,p_1 \in l_1; p,p_2 \in l_2; p',p_1 \in l'_1; p',p_2 \in l'_2;\\
&p,p' \neq p_1,p_2; l_1 \neq l_2; l'_1 \neq l'_2 \}.
\end{align*}
We can pass to the non-degenerate quadrilaterals
$$ 
Q' := \{ (p,p',p_1,p_2, l_1,l_2,l'_1,l'_2) \in Q: p \neq p' \}
$$
since $p = p'$ forces $(p,p_1,p_2,l_1,l_2) = (p',p_1,p_2,l'_1,l'_2)$, 
and so the cardinality of $Q \backslash Q'$ is at most 
$| A| \ll |A|^2/|P|^2 \lesssim | Q|$.  Thus
\be{george}
|Q'| \gtrsim |F|^{4n+4}\ |P|^{-4}.
\end{equation}
Now consider the map $f: Q' \to P^3$ defined by
$$ 
f(p,p',p_1,p_2, l_1,l_2,l'_1,l'_2) := (\frac{1}{2}p + \frac{1}{2}p_1, -p_1 +2 p', \frac{2}{3} p' + \frac{1}{3} p_2).
$$
This triple lies in $l_1 \times l'_1 \times l'_2$ and is thus in 
$P^3$ as claimed.  Now we consider to what extent an element $q \in Q$ 
is determined by $f(q)$.

Write $q = (p,p',p_1,p_2, l_1,l_2,l'_1,l'_2)$, and suppose that $f(q)$ is 
fixed.  From the identity
$$ 
p - p_2 = 2(\frac{1}{2} + \frac{1}{2}p_1) + (-p_1 +2p') - 3 (\frac{2}{3} p' + \frac{1}{3} p_2)
$$
we see that $p-p_2$ is also fixed.  But $p, p_2 \in l_2$, so $p-p_2$ is 
parallel to the direction of $l_2$.  Since the lines in $L$ all point in 
different directions, we thus see that $l_2$ is fixed.  This gives at most 
$|F|$ choices for $p$ (say).  Once $p$ is chosen, the remaining components 
of $q$ are determined since $f(q)$ is fixed (the lines $l_1,l_2,l'_1,l'_2$ 
are then uniquely determined by the distinct points $p,p_1,p_2,p'$).  

To summarize, for each fixed value of $f(q)$ there are at most $|F|$ values 
of $q$.  Since there are at most $|P|^3$ values of $f(q)$, we thus have
$$ |Q'| \lesssim |P|^3\ |F|.$$
Combining this with \eqref{george} we obtain $|P| \lesssim |F|^{(4n+3)/7}$ 
as desired.
\end{proof}

This argument can be combined with Wolff's argument in Proposition \ref{wolff} 
to obtain further improvements, from $(4n+3)/7$ to $(4n+5)/7$ (or better) 
if $\chr(F)$ is sufficiently large.  See \cite{kt:advanced}.  Some of 
these arguments also extend to the maximal function problem, giving 
estimates of a type similar to \eqref{Wolff-est}.  We will not pursue 
these matters here, but remark that anyone interested in improving the 
known Kakeya results in the Euclidean case might first consider the finite 
field case in which several technical difficulties (e.g. small angles and 
small separations) are no longer present\footnote{On the other hand certain 
Euclidean arguments (\cite{katzlaba}, \cite{laba:medium}) require an 
induction on scales type argument which is not directly reproducible for 
finite fields (unless perhaps they have subfields of index 2).}.

\section{The connection between Restriction and Kakeya}\label{conn-sec}

Until now we have pursued the restriction and Kakeya problems independently, 
with only a very slight connection between the two (via Proposition 
\ref{lines}).  In the Euclidean case there is a much stronger relationship 
between the two problems.  For instance, when $S$ is the paraboloid
$$ 
S := \{ (\xi, |\xi|^2): \xi \in \R^{n-1}; |\xi| \lesssim 1 \}
$$
one can use a Taylor approximation at scale $R^{-1/2}$ for some $R \gg 1$ 
to approximate the paraboloid to within $O(1/R)$ by the union of flat disks
$$ 
\bigcup_{\xi \in R^{-1/2} \Z^{n-1}: |\xi| \lesssim 1}
(\xi, |\xi|^2) + \{ (\eta, 2 \xi \cdot \eta): |\eta| \lesssim R^{-1/2} \}.
$$ 
A function supported on one of these disks would have Fourier transform 
concentrated along tubes pointing in the normal direction.  Since all the 
disks point in different directions, we thus begin to see Kakeya-type 
objects appearing.  For more precise discussions we refer to 
\cite{borg:kakeya}, \cite{borg:stein}, \cite{beckner:restrict-endpoint}, 
\cite{vargas:restrict}, \cite{tvv:bilinear}, \cite{tv:cone1}.

To obtain the analogous results for finite fields we cannot work just 
with the standard paraboloid as there is no Taylor approximation in 
finite fields.  Instead, we force the Taylor approximation into existence 
by fiat, by defining the surface
$$ \tilde S := \{ (\xi, \xi \cdot \xi, \eta, \xi \cdot \eta): \eta, \theta \in F^{n-1} \} \subset F^{2n}$$
(we have dropped the factor of $2$ as it makes no difference).

The restriction theory for $S$ and $\tilde S$ are related to each 
other, and to the Kakeya problem, by the following theorem.

\begin{theorem}  Let $2\leq  p,q $ be exponents.  Then we have
\be{r-0}
R^*_S(p \to q) \leq R^*_{\tilde S}(p \to q),
\end{equation}
\be{r-1}
K(\, (q/2)' \to (p/2)'\, ) \leq R^*_{\tilde S}( p \to q )^2\ |F|^{n-1 - \frac{2n}{p}}
\end{equation}
and
\be{r-2}
R^*_{\tilde S}(p \to q) \leq R^*_S(2 \to q)\  K(\, (q/2)' \to (p/2)'\,)^{1/2}.
\end{equation}
\end{theorem}

The estimate \eqref{r-1} shows how restriction estimates imply Kakeya 
estimates, and is the finite field analogue of the result in 
\cite{beckner:restrict-endpoint}.  At the endpoint $p = 2n/(n-1)$
 the result is especially interesting given that there is no loss of 
powers of $|F|$ in this case.  The estimates \eqref{r-2}, \eqref{r-0} 
give a partial converse and is the analogue of the arguments in 
\cite{borg:kakeya}, \cite{borg:stein}.

\begin{proof}  
We first prove \eqref{r-0}.  Let $f$ be any function on $F^n$.  If we 
define the function $\tilde f$ on $F^{2n}$ by
$$ 
\tilde f(x,y) := f(x) \delta_{y,0}
$$
for all $x,y \in F^n$, we see that
$$ 
\widehat{\tilde f}(\eta, \eta \cdot \eta, \theta, \eta \cdot \theta)
= \widehat f(\eta, \eta \cdot \eta)
$$
and hence that
$$ 
\| \widehat{\tilde f}\|_{L^q(\tilde S, d\sigma)} = \| \widehat f\|_{L^q(S, d\sigma)}.
$$
Since $\| \tilde f\|_{L^p(F^{2n}, dxdy)} = \|f\|_{L^p(F^n, dx)}$, the 
claim \eqref{r-0} then follows from \eqref{rest-dual}. 

Now we prove \eqref{r-1}.  By \eqref{kakeya-alt} it suffices to show that
$$
\| \frac 1{|F|^{n-1}} \sum_{v \in F^{n-1}} h(v)\ \chi_{l(x_0(v),v)}\, \|_{L^{p/2}(F^n,dx)}
\leq R^*_{\tilde S}( p \to q )^2\, |F|^{n-1 - \frac{2n}{p}}\,  \| h \|_{L^{q/2}(F^{n-1},dv)}
$$
for all functions $h$ on $F^{n-1}$ and choices of map 
$x_0: F^{n-1} \to F^{n-1}$.

Fix $h$, $x_0$; we may assume that $h$ is non-negative.  Define the 
function $\tilde h$ on $\tilde S$ by
$$ 
\tilde h(\eta, \eta \cdot \eta, \theta,  \eta \cdot \theta) = h(-\eta)^{1/2}\ \ 
e(-x_0(-\eta) \cdot \theta).
$$
Clearly
$$ 
\| \tilde h \|_{L^{q}(\tilde S, d\sigma)} = \| h^{1/2} \|_{L^{q}(F^{n-1}, dv)}
= \| h \|_{L^{q/2}(F^{n-1}, dv)}.
$$
Applying \eqref{rest-dual} we thus obtain
$$
\| (\tilde h d\sigma)\spcheck\|_{L^{p}(F^{2n}, dx)} \leq 
R^*_{\tilde S}( p \to q )\ \|h\|_{L^{q/2}(S, d\sigma)}^{1/2}.$$
For $(\tilde h d\sigma)\spcheck(\underline{x}, x_n, \underline{y}, y_n)$
we get 
$$
 \frac 1{|F|^{2(n-1)}}  \sum_{\eta \in F^{n-1}} \sum_{\theta \in F^{n-1}}
e(\underline{x} \cdot \eta + x_n \eta \cdot \eta + \underline{y} \cdot \theta  + y_n \eta \cdot \theta) 
\ e(-x_0(-\eta) \cdot \theta)\  h(-\eta)^{1/2}.
$$
Performing the $\theta$ summation we get
$$ 
\frac1{|F|^{n-1}} \sum_{\eta \in F^{n-1}} 
e(\underline{x} \cdot \eta + x_n \eta \cdot \eta)\ h(-\eta)^{1/2}\
\delta_{\underline{y},x_0(-\eta)-y_n \eta}.$$
Writing $x = (\underline{x},x_n)$, $y = (\underline{y},y_n)$, this becomes
$$ 
(\tilde h d\sigma)\spcheck(x,y)
= \frac 1{|F|^{n-1}} \sum_{\eta \in F^{n-1}}
e(x \cdot (\eta,\eta \cdot \eta))\ h(-\eta)^{1/2}\
\chi_{l(x_0(-\eta),-\eta)}(y).$$

Making the change of variables $v = -\eta$, we thus have
$$ 
\| \sum_{v \in F^{n-1}} e(x \cdot (-v, v \cdot v))\ h(v)^{1/2}\
\chi_{l(x_0(v),v)}(y) \|_{l^{p}_y l^{p}_x} \leq
R^*_{\tilde S}(p \to q)\ |F|^{n-1}\ \|h\|_{L^{q/2}(S, d\sigma)}^{1/2}.
$$

Since $p \geq 2$, we have from H\"older's inequality that
$$ \| f \|_{l^2_x} \leq |F|^{n(\frac{1}{2}-\frac{1}{p})}\ \| f\|_{l^{p}_x}$$
and hence
$$ 
\| \sum_{v \in F^{n-1}} e(x \cdot (-v,v \cdot v))\ h(v)^{1/2}\
\chi_{l(x_0(v),v)}(y) \|_{l^{p}_y l^2_x} \leq
R^*_{\tilde S}(p \to q)\ |F|^{\frac{n}{2} + \frac{n}{p'} - 1}\ \|h\|_{L^{q/2}(S,d\sigma)}^{1/2}.$$
By orthogonality we have
$$
\| \sum_{v \in F^{n-1}} e(x \cdot (-v,v \cdot v))\ h(v)^{1/2}\
\chi_{l(x_0(v),v)}(y) \|_{l^2_x}
= |F|^{n/2}\ (\sum_{v \in F^{n-1}} h(v)\ \chi_{l(x_0(v),v)}(y)\, )^{1/2}.
$$
Thus we have
$$ 
\| (\sum_{v \in F^{n-1}} h(v)\ \chi_{l(x_0(v),v)}(y)\ )^{1/2} \|_{l^{p}_y} \leq
R^*_{\tilde S}(p \to q)\ |F|^{\frac{n}{p'} - 1}\ \|h\|_{L^{q/2}(S, d\sigma)}^{1/2}.$$
Squaring this and then dividing by $|F|^{n-1}$, we obtain the desired 
estimate \eqref{r-1}.

Finally, we prove \eqref{r-2}, i.e.
$$
\| (h d\sigma)\spcheck\|_{L^{p}(F^{2n}, dx dy)}^2 \leq 
R^*_S(2 \to q)^2\ K(\, (q/2)' \to (p/2)'\,)\ \|h\|_{L^{q}(\tilde S, d\sigma)}^2
$$
for all $h$ on $\tilde S$.

We split $h = \sum_{\alpha \in F^{n-1}} h_\alpha$, where 
$h_\alpha$ is the restriction of $h$ to the ``cap'' 
$$ 
C_\alpha = \{ (\alpha, \alpha \cdot \alpha, \theta, \alpha \cdot \theta): \theta \in F^{n-1} \}.
$$
Fix $y=( \underline{y},y_n)$ and expand
\be{fia}
(h_\alpha d\sigma)\spcheck(x,y) = \frac 1{|F|^{2(n-1)}} \sum_{\theta \in F^{n-1}}
e(x \cdot (\alpha,\alpha \cdot \alpha) + (\underline{y} + \alpha y_n) \cdot \theta)\ \ 
h(\alpha,\alpha \cdot \alpha, \theta, \alpha \cdot \theta).
\end{equation}
We shall write this simply as
\be{fia3} 
(h_\alpha d\sigma)\spcheck(x,y) = \frac 1{|F|^{n-1}}\ e(x \cdot (\alpha,\alpha \cdot \alpha)) \ 
	H(\alpha, \underline {y}+y_n\alpha)\ 
\end{equation}
where $H(\alpha,z)$ is a function which depends on $h$ and $\alpha,
z\in F^{n-1}$, but not on $x$.

Hence by using \eqref{rest-def} we get 
\begin{align*}
\| (h d\sigma)\spcheck\|_{L^{q}(F^{2n},dx dy)}^2 
&= 
\|\  \frac 1{|F|^{n-1}} \sum_{\alpha \in F^{n-1}}  e(x \cdot (\alpha,\alpha \cdot \alpha)) \ 
	H(\alpha, \underline {y}+y_n\alpha)\ \|_{l^q_y l^q_x}^2 \\
&\leq  R^*_S(2 \to q)^2\ \| \ \frac 1{|F|^{n-1}} \sum_{\alpha \in F^{n-1}} 
        |H(\alpha, \underline {y}+y_n\alpha)|^2\  \|_{l^{q/2}_y}.
\end{align*}
For a suitable function $g\in l^{(q/2)'}_y$ with norm $1$ 
we may express the latter squared $l^{q/2}_y$-norm as
$$
\frac 1{|F|^{n-1}} \  
\sum_{\underline{y}\in F^{n-1}}\sum_{y_n\in F} \sum_{\alpha \in F^{n-1}} 
	g(\underline{y},y_n)\ |H(\alpha, \underline {y}+y_n\alpha)|^2.
$$
By replacing $\underline{y}$ by  $\underline{y}-y_n \alpha$ and summing 
over $y_n$ first according to \eqref{max-def} we may estimate this by 
$$
\frac 1{|F|^{n-1}} \  
\sum_{\alpha \in F^{n-1}} g^*(-\alpha)\ \sum_{\underline{y}\in F^{n-1}} |H(\alpha, \underline {y})|^2
$$
which by Plancherel is equal to
$$
\frac 1{|F|^{n-1}} \  
\sum_{\alpha \in F^{n-1}} g^*(-\alpha)\ \frac 1{|F|^{n-1}} \ \sum_{\theta\in F^{n-1}} 
|h(\alpha,\alpha\cdot\alpha,\theta,\alpha\cdot\theta)|^2.
$$
By H\"older's inequality and the Kakeya estimate \eqref{kak-est} we may bound this by
$$
K((q/2)' \to (p/2)')\ \Big( \frac 1{|F|^{n-1}}  
\sum_{\alpha \in F^{n-1}} \big( \frac 1{|F|^{n-1}} \ \sum_{\theta\in F^{n-1}} 
|h(\alpha,\alpha\cdot\alpha,\theta,\alpha\cdot\theta)|^2\  \big)^{p/2}\Big)^{2/p}.
$$
Since $p\geq 2$ we can bound the inner $l^2_\theta$-norm by H\"older's inequality and obtain 
\eqref{r-2} by collecting terms.
\end{proof}

As a sample application of these estimates, we consider the $n=3$ case in 
which $-1$ is not a square.  Interpolating between \eqref{wk} and Proposition 
\ref{cordoba} we have
$$ K(\, 9/4 \to 18/5\,) \lessapprox 1.$$
Combining this with Theorem \ref{BR} and \eqref{r-2} we have
$$ R^*_{\tilde S}(\, 36/13 \to 18/5\, ) \lessapprox 1$$
which by Lemma \ref{ts} implies (since the Fourier transform of 
$d\sigma$ on $\tilde S$ can easily be computed to have some decay)
\be{end}
R^*_{\tilde S}(\, 36/13 \to 18/5+\eps\, ) \leq C_\eps
\end{equation}
for all $\eps > 0$.  

Let us now consider the necessary conditions for $R^*_{\tilde S}(p \to q)$ 
to be bounded.  From \eqref{r-0} the conditions must be at least as 
restrictive as those for $S$, in particular we must have $q \geq 2n/(n-1)$
 from \eqref{q-1}. The condition \eqref{q-2} gives $q \geq np'/(n-1)$, but 
this can be improved to $q \geq (n+1)p'/(n-1)$ by applying \eqref{q-improv} 
with $(n,d,k)$ replaced by $(2n, 2(n-1), n-1)$.  (Note that this latter 
condition is almost obeyed with equality by \eqref{end}).  In analogy with 
the Euclidean restriction conjecture for the paraboloid (which has almost 
the same numerology, except that one must have strict inequality in 
$q > 2n/(n-1)$ because of the multiplicity of scales) one may tentatively 
conjecture that these necessary conditions are also sufficient (at least 
when $-1$ is not a square).  From \eqref{r-1} this conjecture would imply 
the Kakeya set conjecture.  The relevant endpoint is $R^*_{S^*}(2n/(n-1) 
\to 2n/(n-1))$; one can show this endpoint is bounded in two dimensions by 
direct computation but the problem remains open in higher dimensions.  
In three dimensions we see from \eqref{r-2} that this would follow from 
the endpoint restriction bound $R^*_S(2 \to 3)$ for the ordinary 
paraboloid, combined with the endpoint Kakeya bound $K(3 \to 3)$.

\end{document}